\input amstex
\magnification=\magstep1
\documentstyle{amsppt}


\pagewidth{5.375in}
\pageheight{7.443in}

%
%
%


\define\Ftwo{{{\bold F}_{\kern-.125em 2}}}
\define\Fp{{{\bold F}_{\kern-.125em p}}}
\define\Fq{{{\bold F}_{\kern-.125em q}}}
\define\cxcj{{\overline{\phantom{x}}}}
\define\Kb{{\overline{K}}}
\define\kb{{\overline{k}}}
\define\Pb{{\overline{P}}}
\define\pib{{\overline{\pi}}}
\define\Qb{{\overline{Q}}}
\define\Sbar{{\overline{S}}}

\define\Ah{{\widehat{A}}}
\define\Bh{{\widehat{B}}}
\define\Mh{{\widehat{M}}}
\define\Xh{{\widehat{X}}}
\define\Yh{{\widehat{Y}}}

\define\Bap{\pmb{\alpha}_p}
\define\alphatwo{\pmb{\alpha}_2}
\define\BC{{\bold C}}
\define\BQ{{\bold Q}}
\define\BQb{{\overline{\bold Q}}}
\define\BQp{{{\bold Q}_p}}
\define\BR{{\bold R}}
\define\BZ{{\bold Z}}

\define\CB{{\Cal B}}
\define\CC{{\Cal C}}
\define\Krr{{{\Cal K}_{rr}}}
\define\Krl{{{\Cal K}_{r\ell}}}
\define\Klr{{{\Cal K}_{\ell r}}}
\define\Kll{{{\Cal K}_{\ell\ell}}}
\define\CP{{\Cal P}}

\define\Aid{{\frak A}}
\define\aid{{\frak a}}
\define\Bid{{\frak B}}
\define\Oid{{\Cal O}}
\define\pid{{\frak p}}
\define\qid{{\frak q}}
\define\rid{{\frak r}}

\define\charact{\mathop{\roman{char}}\nolimits}
\define\Cl{\mathop{\roman{Cl}}\nolimits}
\define\End{\mathop{\roman{End}}\nolimits}
\define\Gal{\mathop{\roman{Gal}}\nolimits}
\define\Gm{{{\bold G}_m}}
\define\Hom{\mathop{\roman{Hom}}\nolimits}
\define\inv{\mathop{\roman{inv}}\nolimits}
\define\Ker{\mathop{\roman{Ker}}\nolimits}
\define\Mod{\mathop{\roman{Mod}}\nolimits}
\define\Nrd{\mathop{\roman{Nrd}}\nolimits}
\define\opp{{\mathop{\roman{opp}}\nolimits}}
\define\Prin{\mathop{\roman{Pr}}\nolimits}
\define\sep{{\mathop{\roman{sep}}\nolimits}}
\define\Tr{\mathop{\roman{Tr}}\nolimits}
\define\SU{\mathop{\bold{SU}}\nolimits}

\define\Xll{{X_{\ell\ell}}}

\define\rosati{{'}}
\define\eps{\epsilon}

\define\ra{\rightarrow}
\define\lra{\longrightarrow}
\define\llra{\relbar\joinrel\longrightarrow}
\define\lllra{\relbar\joinrel\llra}
\define\llllra{\relbar\joinrel\lllra}
\define\mapright#1{\mathop{\lra}\limits^{#1}}
\define\longmapright#1{\mathop{\llra}\limits^{#1}}
\define\llongmapright#1{\mathop{\lllra}\limits^{#1}}
\define\lllongmapright#1{\mathop{\llllra}\limits^{#1}}

\define\lla{\longleftarrow}
\define\llla{\longleftarrow\joinrel\relbar}
\define\lllla{\llla\joinrel\relbar}
\define\llllla{\lllla\joinrel\relbar}
\define\mapleft#1{\mathop{\lla}\limits^{#1}}
\define\longmapleft#1{\mathop{\llla}\limits^{#1}}
\define\llongmapleft#1{\mathop{\lllla}\limits^{#1}}
\define\lllongmapleft#1{\mathop{\llllla}\limits^{#1}}

\define\Btmapdown#1{\vbox{\vbox to 4pt{}\vbox{\hbox{
             \Big\downarrow\rlap{$\vcenter{\hbox{$\scriptstyle#1$}}$}}\vfill}}}
\define\Btbmapdown#1{\vbox{\vbox to 4pt{}\vbox to 18pt{\hbox{
             \Big\downarrow\rlap{$\vcenter{\hbox{$\scriptstyle#1$}}$}}\vfill}}}

\define\ovl{\overline}

\define\Xgood{{X_{\roman{good}}}}
\define\Ygood{{Y_{\roman{good}}}}
\define\Cgood{{C_{\roman{good}}}}
\define\Cbad{{C_{\roman{bad}}}}
\define\Dgood{{D_{\roman{good}}}}
\define\Dbad{{D_{\roman{bad}}}}

\define\myqed{\hbox to 1em{}\nobreak\hfill$\square$}


\define\SECintro{1}
\define\TTmainPP{1.1}
\define\TTodddim{1.2}
\define\TTmain{1.3}

\define\SECprelim{2}

\define\SECGrot{3}
\define\TTGrotnorms{3.1}
\define\LLlnotp{3.2}
\define\LLp{3.3}
\define\LLpidpart{3.4}
\define\TTGrotiso{3.5}

\define\EEker{(1)}
\define\EEell{(2)}
\define\EEp{(3)}

\define\SECNrd{4}
\define\TTNrd{4.1}
\define\LLKrasner{4.2}
\define\LLLplus{4.3}
\define\LLLsplitsA{4.4}
\define\LLlocalsplit{4.5}
\define\LLlocalglobal{4.6}
\define\TTNrdfirst{4.7}

\define\SECkernels{5}
\define\LLbalanced{5.1}
\define\PPisotropic{5.2}

\define\EEkernels{(4)}

\define\SECgroup{6}
\define\PPreal{6.1}
\define\PPCM{6.2}
\define\CCodddegree{6.3}
\define\PPpushout{6.4}
\define\PPpushouttwo{6.5}
\define\CCfinite{6.6}

\define\EEbottom{(5)}
\define\EEtop{(6)}
\define\EEquot{(7)}
\define\EEfinite{(8)}

\define\SECelement{7}
\define\PPOimage{7.1}
\define\PPzeroconditions{7.2}


\define\refD{1}
\define\refH{2}
\define\refMu{3}
\define\refO{4}
\define\refPR{5}
\define\refRam{6}
\define\refR{7}
\define\refSch{8}
\define\refSe{9}
\define\refT{10}
\define\refW{11}


\define\version{19950827}

 
\message{Version \version}
 

\topmatter
\title
Kernels of polarizations of abelian varieties over finite fields
\endtitle
\rightheadtext{KERNELS OF POLARIZATIONS}
\author
Everett W. Howe
\endauthor
\address
Department of Mathematics, University of Michigan,
Ann Arbor, Michigan 48109\endaddress
\curraddr
Unaffiliated, San Diego, CA 92104\endcurraddr
\email
however\@alumni.caltech.edu
\endemail
\urladdr
http://ewhowe.com
\endurladdr
\dedicatory
Dedicated to Professor Frans Oort
on the occasion of his sixtieth birthday.
\enddedicatory

\thanks
The author was supported in part
by NSA grant number MDA904-95-H-1044
and by a Rackham Faculty Fellowship provided by the 
Office of the Vice-President for Research
and the Horace H.~Rackham Graduate School 
at the University of Michigan.
\endthanks


\abstract
Suppose $\CC$ is an isogeny class of abelian varieties over 
a finite field $k$.  
In this paper we give a partial answer to the question
of which finite group schemes over $k$ 
occur as kernels of polarizations of varieties in $\CC$.
We show that there is an element $I_\CC$ 
of a finite two-torsion group that determines which Jordan-H\"older 
isomorphism classes of finite commutative group schemes over 
$k$ contain kernels of polarizations.  
We indicate how the two-torsion group can be computed from the
characteristic polynomial of the Frobenius endomorphism of the varieties
in $\CC$, and we give some relatively weak sufficient
conditions for the element $I_\CC$ 
to be zero.
Using these conditions, we show that every isogeny class
of simple odd-dimensional abelian varieties over a finite field
contains a principally polarized variety.
As a step in the proofs of these theorems, we prove that if $K$
is a CM-field and $A$ is a central simple 
$K$-algebra with an involution of the second kind, then every 
totally positive real element of $K$ is the reduced norm
of a positive symmetric element of $A$.
\endabstract

\endtopmatter

\document

{\it This is a reproduction of a preprint, dated\/ $27$\,August\/ $1995$, 
of a paper that appeared in the {\rm Journal of Algebraic Geometry} in\/ $1996$.
The only changes between the original\/ $1995$ preprint and this document
are some adjustments to the font size and margins,
the inclusion of this note, the addition of a current address, the updating of
my email address, and the inclusion of a link to my homepage. The numbering of
theorems and so forth is consistent with the published version.} --- EWH, $14$
January $2020$

\head \SECintro. Introduction
\endhead
Every abelian variety over an algebraically closed field is isogenous
to a principally polarized variety, but over non-algebraically-closed
base fields the situation is more complicated.  
In this paper we study the case where the base field is finite,
and we provide some relatively weak sufficient conditions for an
isogeny class of abelian varieties over a finite field to contain
a principally polarized variety.
Our results are easiest to describe for simple varieties.
Suppose $\CC$ is an isogeny class of simple abelian varieties over
a finite field $\Fq$.  The isogeny class $\CC$ is determined by
its {\it Weil number}, which is an algebraic integer $\pi\in\BC$
such that $\pi\pib=q$.
The field $K=\BQ(\pi)$ is either totally real or a so-called
CM-field, that is, a totally imaginary quadratic extension
of a totally real field $K^+$.

\proclaim{(\TTmainPP) Theorem}
 Let $\CC$ and $K$ be as above.
If $K$ is totally real then $\CC$ contains a 
principally polarized variety.  Suppose $K$ is a CM-field.  If
a finite prime of $K^+$ ramifies in $K/K^+$, or if there is a prime of
$K^+$ that divides $\pi-\pib$ and that is inert in $K/K^+$, 
then $\CC$ contains a principally polarized variety.
\endproclaim

{}From Theorem~\TTmainPP\ we obtain an easily-stated result.

\proclaim{(\TTodddim) Theorem}
Every simple odd-dimensional abelian variety over a 
finite field is isogenous to a principally polarized variety.
\endproclaim

The proof of Theorem~\TTmainPP\ involves the idea of an
{\it attainable} group scheme.  If $k$ is a finite field and $\CC$
is an isogeny class of abelian varieties over $k$, we call
a group scheme $X$ over $k$ {\it attainable in $\CC$} if there
is a polarization of a variety in $\CC$ whose kernel is isomorphic
to $X$. The proof of Theorem~\TTmainPP\ depends on our being
able to identify the image of the attainable group schemes in
the Grothendieck group $G(\Ker_\CC)$ of the category $\Ker_\CC$
whose objects are the kernels of isogenies between elements of $\CC$.
To every isogeny class $\CC$ we associate a commutative ring $R$
of a type we will call a {\it real/CM-order};
if $A$ is any variety in $\CC$ then $R$ is isomorphic to the
subring of $\End A$ generated over $\BZ$ by the Frobenius and
Verschiebung endomorphisms of $A$.  In \S\SECGrot\ we
define an isomorphism between $G(\Ker_\CC)$ and the
Grothendieck group $G(\Mod_R)$ of finite-length modules over $R$, 
so that our goal is to identify
the attainable elements of $G(\Mod_R)$, that is, the classes in $G(\Mod_R)$
that contain the images of attainable group schemes. Some simple
considerations show that the attainable elements
lie in a certain subgroup $Z(R)$ of $G(\Mod_R)$ that we define
in \S\SECkernels.  The attainable elements must also be
{\it effective}; that is, an attainable element must be the
class in $G(\Mod_R)$ of an actual $R$-module. In \S\SECkernels\ 
we define a finite two-torsion quotient $\CB(R)$ of $Z(R)$,
and we prove the following theorem.

\proclaim{(\TTmain) Theorem}
There is an element $I_\CC$ of the group $\CB(R)$
such that the elements of $G(\Mod_R)$ that are attainable in $\CC$ are 
precisely the effective elements of $Z(R)$ that map to $I_\CC$ in $\CB(R)$.
In particular, the isogeny class $\CC$ contains a 
principally polarized variety if and only if $I_\CC=0$.
\endproclaim

In order to apply Theorem~\TTmain\ to actual
isogeny classes, we must be able
to calculate the obstruction group $\CB(R)$ and the obstruction element
$I_\CC$.  In \S\SECgroup\ we indicate how $\CB(R)$ can be
calculated from the Weil number of $\CC$. 
In \S\SECelement\ we show that $I_\CC$ must
lie in a certain calculable subgroup of $\CB(R)$.
In many circumstances this subgroup is the zero group (see
Proposition~\PPzeroconditions), and this fact allows us to prove
Theorems~\TTmainPP\ and~\TTodddim.

Stronger results for isogeny classes of ordinary abelian varieties
are obtained in \cite{\refH}. 
The stronger results are made possible by a theorem
of Deligne (\cite{\refD})
that provides an equivalence between the category
of varieties in an ordinary isogeny class
$\CC$ and a certain category of $R$-modules.
When dealing with arbitrary abelian varieties we no longer have
this tool, so one of the main contributions of the present paper
is Theorem~\TTGrotnorms, which says that in general we still at
least have an isomorphism between the Grothendieck groups of $\Ker_\CC$
and $\Mod_R$ that behaves well with respect to endomorphisms.
The other difficulty we face when working with arbitrary varieties
is that their endomorphism algebras are more complicated than those
of ordinary varieties; this is why we must prove the results
on reduced norms found in \S\SECNrd.

\subsubhead Convention
\endsubsubhead
If $A$ and $B$ are varieties over a field $k$, 
then when we speak of a morphism from $A$ to $B$ 
we always mean a $k$-morphism.

\subsubhead
Acknowledgment
\endsubsubhead
The author thanks Vladimir Platonov 
for helpful conversations and suggestions.

\head \SECprelim.  Preliminaries
\endhead
Suppose $k$ is a finite field with $q=p^a$ elements, where $p$ is prime.
If $A$ is a $g$-dimensional abelian variety over $k$,
we let $h_A$ denote the characteristic polynomial of the Frobenius
endomorphism of $A$.  The polynomial $h_A$ is a monic polynomial
of degree $2g$ with integer coefficients, and all of its roots
in the complex numbers have magnitude $q^{1/2}$.
It follows from the theorem of Honda and Tate (see~\cite{\refT}) that
the polynomial $h_A$, the endomorphism algebra $E_A=(\End A)\otimes\BQ$ of
$A$, the center $K_A$ of this algebra, and the subring $R_A$ of $K_A$
that is generated by the Frobenius and Verschiebung 
endomorphisms $F$ and $V$ of $A$,
all depend only on the isogeny class of $A$.  Thus, if $\CC$ is an isogeny
class of abelian varieties over $k$ we can speak of $h_\CC$, $E_\CC$,
$K_\CC$, and
$R_\CC$.  We note that the ring $K_\CC$ is a product of number fields
that is generated over $\BQ$ by the Frobenius endomorphism and
that $E_\CC$ is a product of central simple algebras over
the factors of $K_\CC$.  
By piecing together the reduced norms from the factors of 
$E_\CC$ to the factors of $K_\CC$, we get a homomorphism
from $E_\CC^*$ to $K_\CC^*$ that we denote by $\Nrd$
and continue to call the reduced norm.
Detailed information about the structure of the endomorphism algebras
of abelian varieties over finite fields can be found 
in~\cite{\refW\rm, Ch.~1 and~2}, \cite{\refT},
\cite{\refMu\rm, \S21, Application~I, pp.~193--203}, 
and~\cite{\refRam}.

An {\it order} is a commutative ring that has no nilpotent elements and
that is free and finitely generated as a $\BZ$-module.
Equivalently, an order is a subring of finite index in the integral
closure of $\BZ$ in a finite product of number fields.
An order $R$ is called {\it real} 
if $R\otimes\BQ$ is a product of totally real number fields;
it is called {\it CM} if $R\otimes\BQ$ is a product of CM-fields; and
it is called {\it real/CM} if $R\otimes\BQ$ is a product of real and
CM-fields.
If $R$ is a real/CM-order then there is an involution $\cxcj$ on $R\otimes\BQ$
that is complex conjugation on each factor of $R\otimes\BQ$,
and we say that $R$ is {\it proper} if $\ovl{R}=R$.
We let $R^+$ denote the ring of elements of $R$ that are fixed by 
the involution.
If $\CC$ is an isogeny class of abelian varieties over a finite field,
then the ring $R_\CC=\BZ[F,V]$ is a proper real/CM-order
and the involution $\cxcj$ interchanges $F$ and $V$.

Let $R$ be any commutative ring,
let $M$ be an $R$-module, and let $G$ be an abelian group.
A pairing $e\colon M\times M\ra G$ is called {\it $R$-balanced}
if for all $m$ and $n$ in $M$ and for all $a\in R$ we
have $e(am,n)=e(m,an)$. If $R$ has an involution $\cxcj$,
then the pairing $e$ is called {\it $R$-semi-balanced}
if for all $m$ and $n$ in $M$ and for all $a\in R$ we
have $e(am,n)=e(m,\ovl{a}n)$.

\head \SECGrot. Grothendieck groups
\endhead
Let $p$ be a prime number, let $k$ be a field with $q=p^a$ elements, 
and suppose $\CC$ is an isogeny class of abelian varieties over $k$.
Let $\Ker_\CC$ be the category whose objects are finite commutative 
group schemes that can be embedded (as closed sub-group-schemes) in 
some variety in the isogeny class $\CC$, and whose morphisms are morphisms
of group schemes.  We see that the objects in $\Ker_\CC$ are
those group schemes that can be written $\ker\varphi$ for some
isogeny $\varphi\colon A\ra B$ of elements of $\CC$.
The category $\Ker_\CC$
splits into a product of four subcategories 
whose objects are, respectively, reduced
group schemes whose Cartier duals are reduced, reduced group schemes
whose Cartier duals are local, local group schemes whose Cartier
duals are reduced, and local group schemes whose Cartier duals are
local. We will denote these subcategories by $\Krr$, $\Krl$, $\Klr$,
and $\Kll$, respectively.

Suppose $X$ is an object of $\Krr$ or of $\Krl$.  Then $X$ is completely
determined by the $\Gal(\kb/k)$-module $X(\kb)$, where $\kb$
is an algebraic closure of $k$. The Galois
action on $X(\kb)$ is determined by the action of the $q$th-power
automorphism of $\kb$, which is identical to the action of the Frobenius
endomorphism on $X(\kb)$.  
Since the action of Frobenius on $X$ agrees with the action
of Frobenius on the varieties in $\CC$, the group $X(\kb)$ can be
viewed as a module over the ring $R=R_\CC=\BZ[F,V]$, 
where we let $F$ act as the Frobenius and $V$ as $q/F$.
A morphism between two objects of $\Krr\times\Krl$ corresponds to a 
group homomorphism between their groups of $\kb$-valued points that 
is equivariant under the Frobenius endomorphisms of the two groups.
Thus, the functor from $\Krr\times\Krl$ to the category $\Mod_R$ of
finite $R$-modules that sends a group scheme $X$ to the $R$-module 
$X(\kb)$ is fully faithful and exact.

If $M$ is a finite $R$-module we 
define the {\it dual module\/} of $M$ to be
the $R$-module $\Mh=\Hom_\BZ(M,\BQ/\BZ)$,
where for every $r\in R$ and $\psi\in\Mh$ we let $r\psi\in\Mh$
be defined by $(r\psi)(m)=\psi(\ovl{r}m)$ for all $m\in M$.
Note that if $\pid$ is a maximal ideal of $R$ then $R/\pid$
and $R/\ovl{\pid}$ are dual to one another.
If $X$ is an object of $\Klr$ then its Cartier dual $\Xh$
is an object of $\Krl$, and 
the functor from $\Klr$ to $\Mod_R$ that sends
$X$ to $\widehat{\Xh(\kb)}$ is fully faithful and exact.

Piecing together the functors defined above, we get an exact
fully faithful
functor $\CP$ from the category $\Krr\times\Krl\times\Klr$ to $\Mod_R$.
The functor $\CP$ transforms Cartier duality into the duality of $\Mod_R$
that we defined above, and for every object
$X$ of $\Krr\times\Krl\times\Klr$, the cardinality of
the $R$-module $\CP(X)$ is equal to the rank of the group scheme $X$.

We will need only two facts concerning the category $\Kll$.
The first is that
if it contains any non-zero object then it contains exactly one
simple object, namely the unique simple local-local finite commutative
group scheme $\Bap$.  
The second is that if $\Kll$ contains $\Bap$ then the varieties in
$\CC$ are not ordinary. In this case $F$ and $V$ 
are not coprime in $R$, and there is a unique prime of $R$ that contains
them both: the kernel of the map from $R$ to $\BZ/p\BZ$ that sends both
$F$ and $V$ to zero.

Recall that the {\it Grothendieck group\/} $G(\Mod_R)$ of $\Mod_R$ is
defined to be the quotient of the free abelian group on the
isomorphism classes of objects in $\Mod_R$ by the subgroup generated
by the expressions $M-M'-M''$ for all exact sequences
$0\ra M'\ra M\ra M''\ra 0$ in $\Mod_R$.  We define the Grothendieck
group $G(\Ker_\CC)$ similarly.  
If $M$ is a finite $R$-module we let $[M]_R$ denote its class
in $G(\Mod_R)$; likewise, we let $[X]_\CC$ denote the class 
in $G(\Ker_\CC)$ of an object $X$ of $\Ker_\CC$.
Both $G(\Mod_R)$ and $G(\Ker_\CC)$ are free abelian
groups on the simple objects in their respective categories.
An element of one of these Grothendieck groups is said to be
{\it effective\/} if it is a sum of positive multiples
of simple objects.

Suppose $X$ is a simple object of $\Ker_\CC$. If $X$ is not local-local,
define $\eps([X]_\CC)$ to be the element $[\CP(X)]_R$ of $G(\Mod_R)$;
if $X$ is local-local, then by our previous remarks the group $\BZ/p\BZ$
can be made into an $R$-module $M$ by letting $F$ and $V$ act as zero,
and we define $\eps([X]_\CC)$ to be $[M]_R$.
Extending $\eps$ by linearity, we get a homomorphism from all of
$G(\Ker_\CC)$ to $G(\Mod_R)$.  The exactness of $\CP$ shows
that for all objects $X$ of $\Ker_\CC$ that have no local-local
part, we have $\eps([X]_\CC)=[\CP(X)]_R$.

We know that the product of fields $K=R\otimes\BQ$ is the
localization of $R$ with respect to the multiplicative
set of non-zero-divisors.  
Suppose $\alpha$ is an invertible
element of $K$, say $\alpha=a/b$ where $a$ and $b$ are elements
of $R$ that are not zero-divisors.  The {\it principal element\/}
generated by $\alpha$ is the element
$\Prin_R(\alpha)=[R/aR]_R - [R/bR]_R$ of $G(\Mod_R)$.  It is easy to 
see that $\Prin_R(\alpha)$ is well-defined and that $\Prin_R$ is a 
homomorphism from $K^*$ to $G(\Mod_R)$.

\proclaim{(\TTGrotnorms) Theorem}  Suppose $A$ is a variety in $\CC$,
let $E=(\End A)\otimes\BQ$, and  let $\Nrd\colon E^*\ra K^*$ be the
reduced norm.  If $\alpha$ is an isogeny in $\End A$, then
$$\eps([\ker\alpha]_\CC) = \Prin_R(\Nrd(\alpha)).$$
\endproclaim

\demo{Proof}
The abelian variety $A$ is isogenous to a variety $A'$ 
that may be written
$A'=B_1^{n_1}\cdots B_r^{n_r}$, where the 
$B_i$ are simple abelian varieties
that are not isogenous to one another.  It is not too difficult to
show 
that the theorem will be true for $A$ if and only if it is true
for $A'$, and 
that the theorem will be true for
$A'$ if and only if it is true for each of the varieties $B_i^{n_i}$.
Thus we need only consider the case where $A$ may be written as
the $n$th power of a simple variety $B$.
Since both sides of the theorem's equality are multiplicative, it will be
enough to prove the theorem for a set of endomorphisms
that generate the multiplicative group $E^*$. In particular,
we can restrict our attention to isogenies
$\alpha$ with $\Nrd(\alpha)\in R$. We must show that for every
such $\alpha$ we have
$$\eps([\ker\alpha]_\CC) = [R/\Nrd(\alpha)R]_R. \eqno{\EEker}$$

To show that equality \EEker\ holds we will express both sides
of the equality in terms of the $\BZ$-basis 
$\{[M]_R  : \hbox{$M$ is a simple $R$-module}\}$
of $G(\Mod_R)$. If a particular simple $R$-module
$M$ occurs the same number of times on both sides, we will say
that {\it equality \rom{\EEker} holds at $M$}.  
For a fixed prime number $\ell$,
we will say that {\it equality \rom{\EEker} holds at $\ell$} if
equality  \EEker\ holds at every
simple $R$-module with cardinality a power of $\ell$.
\enddemo

\proclaim{(\LLlnotp) Lemma}
Let $\ell\neq p$ be a prime number.
Then equality \rom{\EEker} holds at~$\ell$.
\endproclaim

\demo{Proof}
Let $h$ be the characteristic polynomial of the Frobenius 
endomorphism of~$B$.  Since $B$ is simple, $h$ is a power of an
irreducible polynomial $P$ of degree $e$, say $h=P^d$.
The endomorphism algebra $E$ of $A$ is a $n\times n$ matrix algebra
over the division algebra $D=(\End B)\otimes \BQ$, which
is of dimension $d^2$ over its center $K$.
If $\pid$ is a prime of $K$ that lies over $\ell$ 
then $D_\pid=D\otimes K_\pid$ is a $d\times d$ matrix algebra 
over~$K_\pid$.  The algebra $E_\pid=E\otimes K_\pid$ is
an $n\times n$ matrix algebra over $D_\pid$, and is thus also an
$nd\times nd$ matrix algebra over $K_\pid$.

Let $T'$ be the $\ell$-adic Tate module of $B$ and let $T$
be the $\ell$-adic Tate module of $A$, so that $T$ is the sum
of $n$ copies of $T'$.
Let $U'=T'\otimes\BQ$ and let $U=T\otimes\BQ$.
The endomorphism algebra $D$ acts on $U'$, so the $\BQ_\ell$-module
$U'$ is also a module over 
$K_\ell=K\otimes\BQ_\ell=\bigoplus_{\pid|\ell}K_\pid$
and $D_\ell=D\otimes\BQ_\ell=\bigoplus_{\pid|\ell}D_\pid$.
As a $K_\ell$-module, $U'$ is free of rank $d$, so that 
$U'=\bigoplus_{\pid|\ell} U'_\pid$ where 
$U'_\pid$ is a $d$-dimensional $K_\pid$ vector space.
The action of the endomorphism algebra $D_\ell$ on $U'$ is 
given by the action of each $d\times d$ matrix algebra
$D_\pid$ on its corresponding vector
space $U'_\pid$.  Similarly, 
$E_\ell=E\otimes\BQ_\ell=\bigoplus_{\pid|\ell} E_\pid$ acts on $U$, and
the action is
given by the action of each $E_\pid$ on its corresponding 
$nd$-dimensional vector space.

Let $X$ be the $\ell$-part of the group scheme $\ker\alpha$ and
let $\alpha_\ell$ be the endomorphism of $U$ obtained from the
endomorphism $\alpha$ of $A$.  Then the group $X(\kb)$
is isomorphic to $T/\alpha_\ell(T)$, and the actions of 
Frobenius and Verschiebung on $X(\kb)$ correspond
to their actions on $T/\alpha_\ell(T)$.  In other words, 
the $R$-module $T/\alpha_\ell(T)$ is isomorphic to $\CP(X)$,
where $\CP$ is the functor we defined above.  Therefore, 
to prove that equality \EEker\ holds at the prime $\ell$  we
need only show that 
$$[T/\alpha_\ell(T)]_{R_\ell}=
      [R_\ell/\Nrd_\ell(\alpha)R_\ell]_{R_\ell},\eqno{\EEell}$$
where $R_\ell=R\otimes\BZ_\ell$, where $\Nrd_\ell$ denotes the
reduced norm from $E_\ell^*$ to $K_\ell^*$, and where we view
$\alpha$ as an element of $E_\ell^*$.

Let $\Oid$ be the ring of integers of $K$, let 
$\Oid_\ell=\Oid\otimes\BZ_\ell$, and let $S$ be any $\Oid_\ell$-lattice
in $U$ that contains $T$.  It is easy to show that 
$[S/\alpha_\ell(S)]_{R_\ell} = [T/\alpha_\ell(T)]_{R_\ell}$
and that 
$[\Oid_\ell/\Nrd_\ell(\alpha)\Oid_\ell]_{R_\ell}=
              [R_\ell/\Nrd_\ell(\alpha)R_\ell]_{R_\ell}$.
Therefore we need only prove equality \EEell\ with $R$
replaced with $\Oid$ and with $T$ replaced with $S$.
Since $\Oid_\ell=\bigoplus_{\pid|\ell} \Oid_\pid$ and since
$S$ splits accordingly into a sum of $\Oid_\pid$-modules
$S_\pid$, we need only prove equality \EEell\ with
$R_\ell$ replaced by $\Oid_\pid$, with $T$ replaced by 
$S_\pid$, and with $\alpha_\ell$ replaced by $\alpha_\pid$,
where $\alpha_\pid$ is the restriction of the endomorphism $\alpha_\ell$
to $U_\pid$.
It follows from \cite{\refSe\rm, \S{}I.5, Lem.~3, p.~17} 
that $[S_\pid/\alpha_\pid(S_\pid)]_{\Oid_\pid}=
 [\Oid_\pid/\det_{\End_{K_\pid}U_\pid}(\alpha_\pid)\Oid_\pid]_{\Oid_\pid}$,
so it will suffice for us to show that 
$\det_{\End_{K_\pid}U_\pid}(\alpha_\pid)=\Nrd_\pid(\alpha)$.
But this last equality follows from the facts that the action of
$E_\pid$ on $U_\pid$ gives us an isomorphism 
$E_\pid\cong\End_{K_\pid}U_\pid$
and that the reduced norm on a matrix algebra is the determinant.
Thus, equality \EEker\ holds at $\ell$.
\myqed\enddemo

\proclaim{(\LLp) Lemma}
Equality \rom{\EEker} holds at $p$.
\endproclaim

\demo{Proof}
Let $D$, $E$, and $h=P^d$ be as in the proof of Lemma~\LLlnotp.
Let $L$ be an unramified complete extension field of $\BQp$ with residue
field $k$, let $\sigma$ be the automorphism of $L$ induced
by the automorphism $x\mapsto x^p$ of $k$,
let $W$ be the ring of integers of $L$, and let $\Aid=W[f,v]$
where $f$ and $v$ are two 
indeterminates subject to the relations $fv=vf=p$, 
$f\lambda=\lambda^\sigma f$,
and $\lambda v=v \lambda^\sigma$ for every $\lambda\in W$.
Let $T'$ be the Dieudonn\'e module
associated to $B$ and $T$ the Dieudonn\'e module associated 
to $A$ (see \cite{\refW\rm, Ch.~1 and~2} or \cite{\refO});
$T'$ and $T$ are left 
$\Aid$-modules, and $T$ is the direct sum of $n$ copies of $T'$.
Let $U'=T'\otimes\BQ$ and $U=T\otimes\BQ$; these
objects are left modules over $\Bid=\Aid\otimes\BQ=L[f,v]$.

Every endomorphism $\beta$ of $B$ gives us an endomorphism $\beta_p$
of $T'$ and of $U'$, so $U'$ is a module over $K_p=K\otimes\BQp$.  
The splitting $K_p=\bigoplus_{\pid|p} K_\pid$ of $K_p$ into a sum
of local fields gives us a splitting $U'=\bigoplus_{\pid|p} U'_\pid$
of $U'$ into  $K_\pid$-modules $U'_\pid$.
The Frobenius endomorphism $F\in K$ of $B$ and the element $f^a$ 
of $\Bid$ (where $a=[k:\Fp]$) have identical actions on $U'$.
If we let $P_\pid$ be the irreducible factor of $P$ in $\BQp[X]$
that is the minimal polynomial over $\BQp$ of $F\in K_\pid$,
then $\Bid$ acts on $U'_\pid$ via reduction to
$\Bid_\pid=\Bid/P_\pid(f^a)\Bid$, which is a central 
simple algebra of dimension $a^2$ over $K_\pid$.  
Similarly, the $\Bid$-module $U$ splits into a direct sum
of submodules $U_\pid$, each of which is a sum of $b$ copies
of $U'_\pid$.

Let $X$ be the $p$-part of the group scheme $\ker\alpha$.
Then the Dieudonn\'e module attached to $X$ is isomorphic
to $T/\alpha_p(T)$.
The simple group schemes occurring in $X$ correspond to
the Jordan-H\"older factors of the left $\Aid$-module $T/\alpha_p(T)$.
Let $T_\pid$ be the image of $T$ in $U_\pid$, and let 
$S=\bigoplus T_\pid \subset \bigoplus U_\pid =U$. It is not hard
to show that the left  $\Aid$-modules $S/\alpha_p(S)$ and $T/\alpha_p(T)$
are Jordan-H\"older isomorphic to one another, and 
$S/\alpha_p(S)\cong\bigoplus T_\pid/\alpha_\pid(T_\pid)$,
where $\alpha_\pid$ is the restriction of $\alpha_p$ to $T_\pid$.
For each $\pid$, let $\Oid_\pid$
denote the ring of integers of $K_\pid$, so that $R_p=R\otimes\BZ_p$ is
a subring of finite index in $\bigoplus\Oid_\pid$; 
it is again not hard to show that 
$[R_p/\Nrd(\alpha)R_p]_R
  = \sum_{\pid|p} [\Oid_\pid/\Nrd_\pid(\alpha)\Oid_\pid]_R$,
where for every term in the sum
we view $\alpha$ as an element of $E_\pid=E\otimes_K K_\pid$ and 
where $\Nrd_\pid$ is the reduced norm from $E_\pid^*$ to $K_\pid^*$.
Therefore, to show that equality \EEker\ holds at $p$ 
we need only prove the following lemma.
\enddemo

\proclaim{(\LLpidpart) Lemma}
Let notation be as above 
and let $\pid$ be any prime of $K$ lying over $p$.
Let $Y$ be the finite group scheme corresponding to the
Dieudonn\'e module $T_\pid/\alpha_\pid(T_\pid)$.  Then
$$\eps([Y]_\CC)=
         [\Oid_\pid/\Nrd_\pid(\alpha)\Oid_\pid]_R.\eqno{\EEp}$$
\endproclaim

\demo{Proof}
Let $\qid$ be the prime of $R$ that lies under the prime $\pid$.
Then the only simple $R$-module that occurs in 
$\Oid_\pid/\Nrd_\pid(\alpha)\Oid_\pid$ is $R/\qid$.
We will show that $R/\qid$ is also the only simple $R$-module
that occurs in $\eps([Y]_\CC)$.

First suppose that $\pid$ does not contain $F$.
Recall that the prime $\pid$ corresponds to an irreducible
factor $P_\pid\in\BQp[X]$ of the polynomial $P$; similarly,
the prime $\qid$ corresponds to an irreducible factor $Q_\qid\in\Fp[X]$
of $P$ modulo $p$, and in fact $(P_\pid \bmod p)$ is a power of $Q_\qid$.
Since $F$ is not in $\pid$ the constant coefficient of $P_\pid$ is
not a multiple of $p$, and the constant coefficient of $Q_\qid$ is
non-zero.
Let $N$ be a simple $\Aid$-module on which $f^a$ satisfies $P_\pid$.
Since $N$ is simple it is killed by $p$, and $v$ must act as zero on $N$
because $fv=p$ and $f$ acts invertibly on $N$.
Also, since a power of $Q_\qid(f^a)$ kills $N$, $f^a$ must in fact
satisfy $Q_\qid$ on $N$.
This shows that $\Aid$ must act on $N$ via reduction to the ring
$$W[f,v]/(p,v,Q_\qid(f^a))=k[f]/(Q_\qid(f^a)),$$
where for every $\lambda\in k$ we have $f\lambda=\lambda^p f$.
This ring is a simple algebra over its center, which is
$\Fp[f^a]/(Q_\qid(f^a))\cong \Fp[X]/(Q_\qid(X))\cong R/\qid$.
There is exactly one simple representation of this algebra,
so $N$ is the {\it unique} simple left $\Aid$-module on which
$f^a$ satisfies $P_\pid$.  Let $Z$ be the simple reduced-local
group scheme corresponding to $N$.  One can show that the
action of the Frobenius on $Z(\kb)$ satisfies the polynomial $Q_\qid$,
so $\CP(Z)\cong R/\qid$.  Thus, $R/\qid$ is the only
simple $R$-module occurring in $\eps([Y]_\CC)$.

Next, suppose $\pid$ does not contain $V$, so that the conjugate
ideal $\ovl{\pid}$ of $\pid$ does not contain $F$.
If we take duals everywhere
and apply the preceding argument, we find that $R/\ovl{\qid}$ is
the only simple $R$-module that occurs in $\eps([\Yh]_\CC)$.
Taking duals again, we find that $R/\qid$ is the only
simple $R$-module occurring in $\eps([Y]_\CC)$.

Finally, suppose $\pid$ contains both $F$ and $V$, 
so that $\qid$ is the unique prime of $R$ that contains $F$ and $V$
and $R/\qid\cong \BZ/p\BZ$.
Let $N$ be any simple quotient of the $\Aid$-module $T_\pid$.  
Then $f$ and $v$ must
both act as zero on $N$, so $N$ must be the simple $\Aid$-module $k$.
This $\Aid$-module corresponds to the local-local group scheme $\Bap$,
and we know $\eps([\Bap]_\CC)=[R/\qid]_R$;
thus, $R/\qid$ is the only
simple $R$-module occurring in $\eps([Y]_\CC)$.

In every case $R/\qid$ is the only $R$-module occurring
on either side of \EEp. Let $N$ be the unique
simple left $\Aid$-module that is a quotient of $T_\pid$.
We must show that the number of occurrences of $N$ in
$T_\pid/\alpha_\pid(T_\pid)$ is equal that 
of $R/\qid$ in $\Oid_\pid/\Nrd_\pid(\alpha)\Oid_\pid$.
We can find these numbers simply by looking at the cardinalities
of the objects involved.

Let $m$ denote the $W$-length of $N$. 
Since the unique simple $W$-module is $k$ and $\#k=p^a$, the cardinality
of $N$ is  $p^{am}$. 
On the other hand, the cardinality of $R/\qid$ is equal to the 
rank of the group scheme corresponding to $N$, and the
general Dieudonn\'e theory tells us that this rank is $p^m$; thus
we have $\#R/\qid=p^m$. Therefore, to finish
the proof of the lemma we must show that the cardinality
of $T_\pid/\alpha_\pid(T_\pid)$ is the $a$th power of that of
$\Oid_\pid/\Nrd_\pid(\alpha)\Oid_\pid$.

To aid us in computing $T_\pid/\alpha_\pid(T_\pid)$ we choose
an $\Oid_\pid$-lattice $S$ in $U_\pid$.  It is not hard to
show that $\#T/\alpha_\pid(T)=\#S/\alpha_\pid(S)$, and a standard
lemma shows that the $\Oid_\pid$-modules $S/\alpha_\pid(S)$ and
$\Oid_\pid/\det_{\End_{K_\pid}U_\pid}(\alpha_\pid)\Oid_\pid$ 
are Jordan-H\"older
isomorphic; here we view $U_\pid$ as a $K_\pid$-vector space and we view
$\alpha_\pid$ as a $K_\pid$-linear automorphism of $U_\pid$.
Thus to prove the lemma
we must show that the cardinality of 
$\Oid_\pid/ \det_{\End_{K_\pid}U_\pid}(\alpha_\pid)\Oid_\pid$
is the $a$th power of that of $\Oid_\pid/\Nrd_\pid(\alpha)\Oid_\pid$,
and to prove this we need only show that
$\det_{\End_{K_\pid}U_\pid}(\alpha_\pid)=\Nrd_\pid(\alpha)^a$.

We know from \cite{\refW\rm, p.~527} that the dimension of $U'_\pid$ as
a $K_\pid$-vector space is $ad$.  Since $U_\pid$ is the sum of $n$
copies of $U'_\pid$, we have $\dim_{K_\pid} U_\pid = and$.  
On the other hand, $D_\pid$ is a
central simple algebra of dimension $d^2$ over $K_\pid$, so
$E_\pid$ is a central simple algebra of dimension $(nd)^2$ 
over $K_\pid$.  Let $\Kb_\pid$ be an algebraic closure of $K_\pid$.
The $\Kb_\pid$-algebra $E_\pid\otimes_{K_\pid} \Kb_\pid$ is an 
$nd\times nd$ matrix algebra, the map 
$E_\pid\otimes_{K_\pid} \Kb_\pid\ra
  \End_{\Kb_\pid}(U_\pid\otimes_{K_\pid} \Kb_\pid)$ is an
embedding of a $nd\times nd$ matrix algebra into an $and\times and$ matrix
algebra, and the reduced norm on $E_\pid\otimes_{K_\pid}\Kb_\pid$
is the determinant map.
The embedding of the first algebra into the second differs from the
block diagonal embedding by an inner automorphism.  It is clear
that if we map an element $\beta$ of $E_\pid\otimes_{K_\pid}\Kb_\pid$
to the second algebra and 
take its determinant, we get the $a$th power of the
determinant of $\beta$ as an element of $E_\pid\otimes_{K_\pid}\Kb_\pid$.
Thus $\det_{\End_{K_\pid}U_\pid}(\alpha_\pid)=\Nrd_\pid(\alpha)^a$,
and the proofs of Lemmas~\LLpidpart\ and~\LLp\ are complete.
\myqed\enddemo

Together, Lemmas~\LLlnotp\ and ~\LLp\ complete the proof of
Theorem~\TTGrotnorms.\myqed

\proclaim{(\TTGrotiso) Theorem}
The homomorphism $\eps$ is an isomorphism.
\endproclaim

\demo{Proof} 
We know that $\eps$ takes the class of a simple group
scheme to the class of a simple $R$-module, so if $\eps$ were not an
isomorphism there would be a simple $R$-module $M$ such that
$M$ does not occur in $\eps(P)$ for any $P\in G(\Ker_\CC)$. But if
$M$ is a simple $R$-module killed by a prime $\ell$ then $M$
occurs in $\eps(\ker\ell)$ by Theorem~\TTGrotnorms.  Thus $\eps$ is
an isomorphism. 
\myqed\enddemo

\head \SECNrd. Reduced norms
\endhead    
Suppose $K$ is a CM-field whose maximal
real subfield $K^+$ has degree $e_0$ over $\BQ$.  
Let $\sigma$ be the non-trivial automorphism
of $K/K^+$.  Let $A$ be a central simple algebra over $K$ and let
$\tau$ be an involution of $A$ of the second kind --- that is,
$\tau$ should extend the involution $\sigma$ of $K$. 
An element $x$ of $A$ is {\it symmetric} if $\tau(x)=x$, and we
denote by $A^+$ the set of symmetric elements of $A$.
One can show that all of the roots
of the characteristic polynomial 
of a symmetric $x\in A$ are totally real algebraic numbers.
We call a symmetric $x\in A$ {\it totally positive} if all of the
roots of its characteristic polynomial are positive 
under every embedding of $\BQb$ into $\BC$.
An $x$ in $A^+$ is totally positive if and only if 
$x$ is sent to a positive real number under every 
ring homomorphism from $K^+(x)$ to $\BC$.  
We denote the sets of totally positive elements of $A^+$ and of $K^+$ 
by $A_0$ and $K_0$, respectively.  The 
Hasse-Schilling-Maass theorem (\cite{\refR\rm, Thm.~33.15, p.~289})
shows that $\Nrd(A^*)=K^*$, where $\Nrd$ is
 the reduced norm from $A^*$ to $K^*$.  
Our next result, which is related to 
\cite{\refPR\rm, \S6.7, Ex.~2, p.~368}, describes the set
of reduced norms of totally positive elements.

\proclaim{(\TTNrd) Theorem}
Let $K$ be a CM-field and let $A$ be a central
simple $K$-algebra with an involution $\tau$ of the second kind
that fixes $K^+$.  Then $\Nrd(A_0)=K_0$.
\endproclaim

\demo{Proof}
We note for future reference that the fact that $A$ has an involution
of the second kind gives us information about its local invariants.
Landherr's theorem (\cite{\refSch\rm, Thm.~10.2.4, p.~355})
tells us that if $\qid$ is a prime of $K$ such that 
$\qid=\ovl{\qid}$ then $\inv_\qid(A)=0$, while if $\qid\neq\ovl{\qid}$
then $\inv_\qid(A)+\inv_{\ovl{\qid}}(A)=0$. 

It is clear that $\Nrd(A_0)\subset K_0$, so we are left to show
that for every $y\in K_0$ there is an $x\in A_0$ with $\Nrd(x)=y$.  
The first step of our proof will be to find a monic polynomial
of degree $d=\sqrt{[A:K]}$ that has constant term $(-1)^dy$ and
that satisfies local conditions at a set $S$ of primes of $K^+$
that we now define.
The polynomial we construct will depend on $A$ only insofar
as $A$ determines $S$.
Let $S_\infty$ be the set of infinite primes of $K^+$.
Let $S_1$ be the set of finite primes $\pid$ of $K^+$ that lie under $\qid$
of $K$ with $\inv_\qid(A)\neq0$; Landherr's theorem tells us that
every such $\pid$ splits in $K/K^+$.
If the set $S_1$ is empty, replace it with a set consisting of one finite
prime of $K^+$ that splits in $K/K^+$.
Let $S_2$ be the set of finite primes $\pid$ of $K^+$ such that $\pid$
does not split in $K/K^+$ and such that the algebraic group
$\SU(A,\tau)$ over $K^+$ is not $K^+_\pid$-quasisplit
(that is, no geometric Borel subgroup of 
$\SU(A,\tau)\otimes_{K^+}K^+_\pid$
should be definable over $K^+_\pid$).
The set $S_2$  is finite by \cite{\refPR\rm, Thm.~6.7, p.~291}.
For each $\pid\in S=S_\infty\cup S_1\cup S_2$ we will define 
a monic polynomial $f_\pid$ in $K^+_\pid[X]$ of degree~$d$.

Suppose $\pid\in S_\infty$.  Let $\varphi\colon K^+\ra K^+_\pid=\BR$ be the 
embedding of $K^+$ into $\BR$ corresponding to $\pid$. By assumption
$\varphi(y)>0$.  Let $\alpha_1,\ldots,\alpha_d$
be $d$ distinct positive real numbers whose product is $y$,
and let $f_\pid=(X-\alpha_1)\cdots(X-\alpha_d)$.
Clearly $f_\pid$ has constant coefficient $(-1)^d y$.

Suppose $\pid\in S_1$.  Let $f_\pid$ be any monic irreducible polynomial
in $K^+_\pid[X]$ with degree $d$ and constant coefficient $(-1)^d y$;
such polynomials exist by \cite{\refR\rm, Cor.~33.13, p.~288}.
 
Suppose $\pid\in S_2$.  Let $\alpha_1,\ldots,\alpha_d$ be distinct
elements of $K^+_\pid$ whose product is $y$ and let 
$f_\pid=(X-\alpha_1)\cdots(X-\alpha_d)$.
Clearly the constant term of $f_\pid$ is $(-1)^d y$.

We want to find a polynomial $f\in K^+[X]$ that approximates
each of the $f_\pid$ that we have just defined.
The next lemma, which is a slight generalization of
Krasner's lemma (see \cite{\refR\rm, Lem.~33.8, p.~284}),
will give us a neighborhood of the coefficients of $f_\pid$ for
all $\pid\in S$.

\proclaim{(\LLKrasner) Lemma}
Let $K$ be a field complete with respect to a 
valuation and let $f=X^n+a_{n-1}X^{n-1}+\cdots+a_0$
be a separable monic polynomial in $K[X]$.
Let $L$ be the $n$-dimensional $K$-algebra $K[X]/(f)$,
let $\alpha$ be the image of $X$ in $L$, 
and use the valuation on $K$ to define a topology on $L$.  
Let $U$ be a neighborhood of $\alpha$ in $L$.  Then there is
a neighborhood $V$ of $(a_0, a_1, \ldots, a_{n-1})\in K^n$ such that
for any $(b_0, b_1, \ldots, b_{n-1})$ in $V$ the following statement
is true:
Let $g=X^n+b_{n-1}X^{n-1}+\cdots+b_0$, let $M=K[X]/(g)$, and let $\beta$
be the image of $X$ in $M$.  Then there is a $K$-algebra isomorphism 
$i\colon M\ra L$ such that $i(\beta)\in U$.
\endproclaim

\demo{Proof}
Let $c\colon L\ra K^n$ be the continuous map that sends an element
of $L$ to the $n$-tuple of the coefficients of its
characteristic polynomial
over $K$, and let $\ovl{c}\colon L\otimes_K K^\sep \ra (K^\sep)^n$ 
be the map obtained from $c$, where $K^\sep$ is a
separable closure of $K$.  
Applying \cite{\refPR\rm, Lem.~6.25, p.~363)},
we see that $\ovl{c}$ is a local isomorphism near $\alpha$,
and we can find a $\Gal(K^\sep/K)$-invariant neighborhood 
$W\subset L\otimes_K K^\sep$ of $\alpha$ on which $\ovl{c}$ is
an isomorphism.  Thus, we can take $V$ to be $c(U\cap W)$.
\myqed\enddemo

Suppose $\pid$ is a prime in $S$.  Apply Lemma~\LLKrasner\
to the field $K_\pid^+$ and the polynomial $f_\pid$.
If $\pid$ is a finite prime, 
take $U$ to be the open set $L$ (in the notation of the lemma),
while if $\pid$ is infinite take $U$ to be the
open subset $(0,\infty)^d$ of $L\cong \BR^d$.
Let $V_\pid$ be the neighborhood of the coefficients of
$f_\pid$ that is given by the lemma.
Let $f$ be any monic degree $d$ polynomial in $K^+[X]$ that has 
constant term $(-1)^d y$ and whose coefficients lie in 
the neighborhood $V_\pid$ for every $\pid\in S$.  Let $L^+=K^+[X]/(f)$,
and let $x$ be the image of $X$ in $L^+$.

\proclaim{(\LLLplus) Lemma}
The $K^+$-algebra $L^+$ is a totally real number field,
$x$ is totally positive, and $N_{L^+/K^+}(x)=y$.
\endproclaim

\demo{Proof}  The polynomial $f$ is irreducible, because for every
$\pid$ in the non-empty set $S_1$ we have 
$K^+_\pid[X]/(f)\cong K^+_\pid[X]/(f_\pid)$, and $f_\pid$ is irreducible.
Therefore $L^+$ is a field.  The field $L^+$ is totally real because for every
infinite prime $\pid$ of $K^+$ we have an isomorphism
$i_\pid\colon K^+_\pid[X]/(f)\cong K^+_\pid[X]/(f_\pid)\cong \BR^d$.
The element $x$ of $L^+$ is totally positive because for every
infinite $\pid$ Lemma~\LLKrasner\ says
that $i_\pid(x)$ lies in the subset $(0,\infty)^d$ of $\BR^d$.
And finally, the constant term of $f$ tells us that $N_{L^+/K^+}(x)=y$.
\myqed\enddemo

The compositum of the fields $L^+$ and $K$ is a CM-field $L$
with maximal real subfield $L^+$, and the non-trivial automorphism
$\rho$ of $L$ over $L^+$ extends the automorphism $\sigma$ of $K$.
Clearly $x$ is fixed by $\rho$.
We will show in several steps that we can embed $(L,\rho)$
in $(A,\tau)$ as a $K$-algebra with involution.  

\proclaim{(\LLLsplitsA) Lemma}
There is an embedding $L\ra A$
of $K$-algebras without involution.
\endproclaim

\demo{Proof} 
Since $L$ has degree $d$ over $K$ and $\dim_K(A)=d^2$, we see by
\cite{\refR\rm, Cor.~28.10, p.~240} that we can embed $L$ in $A$ if $L$ is a 
splitting field for $A$.  By \cite{\refR\rm, Thm.~32.15, p.~278}, the field
$L$ splits $A$ if and only if for every prime $\rid$ of $L$ we
have $(L_\rid:K_\qid)\inv_\qid(A) = 0$ in the Brauer group of $K_\qid$,
where $\qid$ is the prime of $K$ lying under $\rid$.  Since $\inv_\qid(A)$
is killed by $d$, it will be enough for us to show that if $\qid$ is a 
prime of $K$ with $\inv_\qid(A)\neq 0$ then $\qid$ is inert in $L$,
for then there will be only one prime $\rid$ of $L$ over $\qid$
and $(L_\rid:K_\qid)=d$.

So suppose $\qid$ is a prime of $K$ with $\inv_\qid(A)\neq 0$.
Let $\pid$ be the prime of $K^+$ lying under $\qid$; the prime $\pid$
is an element of the set $S_1$.
By construction we have $K^+_\pid[X]/(f)\cong K^+_\pid[X]/(f_\pid)$,
and since $f_\pid$ is irreducible for $\pid\in S_1$ this last ring
is a field; in other words, the prime $\pid$ of $K^+$ is inert in
$L^+=K^+[X]/(f)$.  But $\pid$ splits in $K/K^+$, so $\qid$ is inert
in $L/K$. 
\myqed\enddemo

\proclaim{(\LLlocalsplit) Lemma}
For every prime $\pid$ of $K^+$,
there is an embedding 
$$(L\otimes_{K^+} K^+_\pid, \rho)\ra (A\otimes_{K^+} K^+_\pid,\tau)$$
of $K^+_\pid$-algebras with involution.
\endproclaim

\demo{Proof}
First suppose $\pid$ is an infinite prime, and let $\qid$ be the
prime of $K$ lying over $\pid$. Then
there is an isomorphism 
$L\otimes_{K^+}K^+_\pid \ra \BC^d$
that takes $\rho$ to the involution that is complex conjugation on
each copy of $\BC$.  There is also an isomorphism 
$A\otimes_{K^+}K^+_\pid\ra M_d(\BC)$ 
that takes $\tau$ to 
the conjugate transpose involution.  Then the first of these algebras
with involution is isomorphic to the diagonal of the second.

Next suppose $\pid$ is a finite prime that splits into distinct primes
$\qid$ and $\ovl{\qid}$ of $K$.  
Then there is an isomorphism
$i\colon A\otimes_{K^+}K^+_\pid\ra 
          (A\otimes_K K_\qid)\times(A\otimes_K K_{\ovl{\qid}})$.
Let $B$ be the $K^+_\pid$-algebra $A\otimes_K K_\qid$.  
Since $\inv_\qid(A)+\inv_{\ovl{\qid}}(A)=0$, the $K^+_\pid$-algebra
$A\otimes_K K_{\ovl{\qid}}$ is isomorphic to the opposite algebra
of $B$, and the isomorphism $i$ takes $\tau$ to the involution of
$B\times B^\opp$ that switches factors.
Similarly, there is an isomorphism 
$L\otimes_{K^+}K^+_\pid \ra 
          (L\otimes_K K_\qid)\times(L\otimes_K K_{\ovl{\qid}})$
that takes $\rho$ to the involution that switches factors.
Therefore, to show that $(L\otimes_{K^+}K^+_\pid , \rho)$
embeds in $(A\otimes_{K^+}K^+_\pid, \tau)$
we need only show that $L\otimes_K K_\qid$ embeds in $A\otimes_K K_\qid$.
But we already know that $L$ embeds in $A$, so we are done.

Next suppose $\pid$ is a finite prime at which $\SU(A,\tau)$ 
is quasisplit.
Then the proof of \cite{\refPR\rm, Lem.~6.26, p.~365} shows that
$(L\otimes_{K^+} K^+_\pid, \rho)$
can be embedded in $(A\otimes_{K^+} K^+_\pid,\tau)$.

We are left to consider the finite primes of $K^+$ that
do not split in $K$ and at which $\SU(A,\tau)$ is not quasisplit.
These are exactly the primes in $S_2$.  
Let $\pid$ be one such prime and let $\qid$ be the prime of $K$ 
lying over $\pid$.
Now, $(L\otimes_{K^+} K^+_\pid, \rho)$ is isomorphic to
$(K_\qid[X]/(f), \sigma)$, where $\sigma$ acts trivially on $X$,
and this last algebra with involution is isomorphic to 
$(K_\qid[X]/(f_\pid), \sigma)$.  Thus, to embed 
$(L\otimes_{K^+}K^+_\pid, \rho)$ into $(A\otimes_{K^+}K^+_\pid, \tau)$
we must find an element $\alpha$ of $A\otimes_{K^+}K^+_\pid$ that
is fixed by $\tau$ and that has minimal polynomial~$f_\pid$.

By Landherr's theorem $\inv_\qid(A)=0$, so 
$A\otimes_{K^+}K^+_\pid=A\otimes_K K_\qid$ is a matrix algebra over
$K_\qid$.  From \cite{\refSch\rm, Thm.~8.7.4, pp.~301--302}
and \cite{\refSch\rm, Thm.~7.6.3, p.~259} 
we see that there is an isomorphism $i\colon A\otimes_K K_\qid\ra M_d(K_\qid)$ 
and a diagonal matrix $a\in M_d(K_\qid)$ with entries
in $K^+_\pid$ such that $i$ takes the
involution $\tau$ to the involution $\eta$ of $M_d(K_\qid)$ defined by
$\eta(x)=a x^* a^{-1}$, where $x^*$ is the conjugate transpose of $x$.
Recall that $f_\pid$ was defined to be $(X-\alpha_1)\cdots (X-\alpha_d)$
where the $\alpha_i$ were distinct elements of $K^+_\pid$.
Let $\alpha$ be the diagonal matrix 
$\langle \alpha_1,\ldots,\alpha_d \rangle\in M_d(K_\qid)$.
Then $\alpha$ is fixed by $\eta$ and
has minimal polynomial $f_\pid$.  Using this $\alpha$,
we can embed $(L\otimes_{K^+} K^+_\pid, \rho)$ into 
$(A\otimes_{K^+}K^+_\pid, \tau)$, and we are done.
\myqed\enddemo

The notation of the next lemma is independent of that of the
rest of this section.
The notation is chosen to agree with that of the lemma's proof, which is to be
found in~\cite{\refPR}.

\proclaim{(\LLlocalglobal) Lemma}
Suppose $L/K$ is a quadratic extension
of number fields, suppose $F$ is an extension of $K$ that is 
linearly disjoint from $L$ over $K$, suppose $A$ is a central
simple $L$-algebra with an involution $\tau$
of the second kind that fixes $K$, and suppose $(F:K)^2=(A:L)$.
Let $P$ be the compositum of $F$ and $L$, so that the non-trivial
automorphism of $L$ over $K$ extends to an involution $\sigma$ of $P$ that 
fixes $F$.  Then there is an embedding $\theta\colon(P,\sigma)\ra(A,\tau)$
of algebras with involution if and only if there is an embedding
$\epsilon\colon P\ra A$ of algebras without involution and for each prime
$\pid$ of $K$ there exists an embedding 
$\theta_\pid\colon(P\otimes_K K_\pid, \sigma)\ra(A\otimes_K K_\pid, \tau)$
of algebras with involution.
\endproclaim

\demo{Proof}  The `only if' statement is clear.
The `if' statement is proven in~\cite{\refPR}. The argument is
found at the end of \S6.7, from the last paragraph of page~366
onward. 
\myqed\enddemo

Together, Lemmas~\LLLsplitsA,~\LLlocalsplit, and~\LLlocalglobal\ show
that there is an embedding $(L,\rho)\ra(A,\tau)$ of algebras with
involution.  View $x\in L$ as an element of $A$ via this embedding.
Then $x$ is fixed by $\tau$, and in fact
$x$ is a totally positive element of $A$. Furthermore,
$\Nrd(x)=N_{L/K}(x)=N_{L^+/K^+}(x)=y$ by Lemma~\LLLplus,
so Theorem~\TTNrd\ is proved.
\myqed\enddemo

The second theorem we will prove in this section
involves central simple algebras over a totally real number field.
Let $K$ be a totally real number field of degree $e$ over $\BQ$
and let $A$ be a central
simple $K$-algebra with an involution $\tau$ of the first kind ---
that is, $\tau$ should be the identity on $K$.  
Let $A^+$ denote the set of elements of $A$ that are fixed by $\tau$.
Suppose that the form  $x\mapsto\Tr_{A/\BQ}(x\tau(x))$ is 
positive definite,
and suppose that $A$ is ramified at all of the infinite primes of $K$.  
Then $A$ is isomorphic
to a matrix algebra over a quaternion algebra $D$ over $K$,
say $A\cong M_d(D)$,
and one can use this fact to show that all of the roots of the
characteristic polynomial of a symmetric element of $A$ are totally
real algebraic numbers.
We call an element $x$ of $A^+$ 
{\it totally positive} if all of the
roots of its characteristic polynomial are positive 
under every embedding of $\BQb$ into $\BC$.
An $x$ in $A^+$ is totally positive if and only if 
$x$ is sent to a positive real number under every 
ring homomorphism from $K(x)$ to $\BC$.  
We denote the sets of totally
positive elements of $A^+$ and of $K$ by $A_0$ and $K_0$, respectively.
The Hasse-Schilling-Maass theorem shows that $\Nrd(A^*)=K_0$. 

\proclaim{(\TTNrdfirst) Theorem}
Let $K$ be a totally real number field
and let $A$ be a central simple $K$-algebra with an involution $\tau$
of the first kind such that $x\mapsto\Tr_{A/\BQ}(x\tau(x))$ is
positive definite.  Suppose $A$ is ramified at all of the infinite
primes of $K$.  Then $\Nrd(A_0)=(K_0)^2$.
\endproclaim

\demo{Proof} 
Let $x\mapsto x^*$ be the conjugate transpose involution on $M_d(D)$.
By \cite{\refSch\rm, Thm.~8.7.4, pp.~301--302} and 
\cite{\refSch\rm, Thm.~7.6.3, p.~259}, there is an isomorphism
$i\colon A\ra M_d(D)$ and a diagonal matrix $a\in M_d(D)$ with $a^*=\pm a$
such that the isomorphism $i$ takes the
involution $\tau$ to the involution $\eta$ of $M_d(D)$ defined by
$\eta(x)=a x^* a^{-1}$.
The argument at the bottom of page~195 of~\cite{\refMu} show that $a^*=-a$
is impossible when $A$ ramifies at an infinite prime,
so we must have $a^*=a$. This shows that $a$ is a diagonal matrix 
with entries in $K$.

It is now easy to show that $\Nrd(A_0)\supset(K_0)^2$.
Suppose $y\in K_0$.  Let $x\in M_d(D)$ be the diagonal matrix with
$y$ in the upper left corner and with ones elsewhere on the diagonal.
Clearly $i^{-1}(x)$ is fixed by $\tau$ and is totally positive,
and $\Nrd(i^{-1}(x))=\Nrd_{D/K}(y)=y^2$.
Therefore $\Nrd(A_0)\supset(K_0)^2$.

Now we show the reverse inclusion holds.  Suppose $x\in A$ is 
fixed by $\tau$ and is totally
positive. By \cite{\refSch\rm, Thm.~7.6.3, p.~259}, there is an
invertible element $b$ of $M_d(D)$ such that $c=bi(x)\eta(b)$ is 
diagonal, and since $i(x)$ and $c$ represent isomorphic Hermitian forms
on a $d$-dimensional $D$-vector space, $c$ must be totally positive.
Write the diagonal matrix $c$ as $\langle c_1,\ldots,c_d\rangle$
with $c_i\in D$.
Since $c$ is fixed by $\eta$, each of the elements $c_i$ must be fixed
by the standard involution of $D$, so each $c_i$ is in $K$.
Furthermore, since $c$ is totally positive, each $c_i$ must lie in $K_0$.
Thus, the reduced norm of $c$, which is $(c_1\cdots c_d)^2$, is in 
$(K_0)^2$.  Since 
$$\eqalign{\Nrd_{M_d(D)/K}(c)&=\Nrd(x)\Nrd_{M_d(D)/K}(b)
                                      \Nrd_{M_d(D)/K}(\eta(b))\cr
                             &=\Nrd(x)(\Nrd_{M_d(D)/K}(b))^2}$$
and since $\Nrd_{M_d(D)/K}(b)\in K_0$ by Hasse-Schilling-Maass,
we see that $\Nrd(x)\in(K_0)^2$. 
\myqed\enddemo

\head \SECkernels. Kernels of polarizations
\endhead
We begin this section by defining the obstruction group $\CB(R)$.
Suppose $R$ is a proper real/CM-order. 
In \S\SECGrot\ we defined a duality $M\mapsto \Mh=\Hom_\BZ(M,\BQ/\BZ)$ 
on $\Mod_R$.  We use this duality to define an involution $\cxcj$
on $G(\Mod_R)$ by setting $\ovl{[M]_R}=[\Mh]_R$ for $R$-modules $M$.
An element $P$ of $G(\Mod_R)$ is {\it symmetric\/} if $P=\Pb$. 
By considering finite $R$-modules to be finite $R^+$-modules,
we obtain the {\it norm} homomorphism $N_{R/R^+}$ from
$G(\Mod_R)$ to $G(\Mod_{R^+})$, 
which satisfies $N_{R/R^+}([M]_R)=[M]_{R^+}$
for every finite $R$-module $M$.  We define $Z(R)$ to be the 
set of symmetric elements of the kernel of the homomorphism
$G(\Mod_R)\ra G(\Mod_{R^+})\otimes(\BZ/2\BZ)$ obtained from the norm.

Let $K=R\otimes\BQ$ and write $K=K_1\times\cdots\times K_r$
as a product of number fields.  If $K_i$ is a CM-field define $K_i^\dagger$
to be the multiplicative group of totally positive elements of $K_i^+$, and
if $K_i$ is real let $K_i^\dagger$ be the group of squares of totally
positive elements of $K_i$. Let $K^\dagger$ be the group
$K_1^\dagger\times\cdots\times K_r^\dagger$.  Consider the subgroup
$\Prin_R(K^\dagger)$ of $G(\Mod_R)$,
where $\Prin_R$ is the homomorphism that was defined in \S\SECGrot.
Clearly the elements of $\Prin_R(K^\dagger)$ are symmetric, and
in fact $\Prin_R(K^\dagger)\subset Z(R)$.
It suffices to verify this in the case where $R$ is a domain.
In this case we must show that for every $a\in R^+\cap K^\dagger$,
the image of $R/aR$ in $G(\Mod_{R^+})$ lies in $2G(\Mod_{R^+})$.
If $R$ is real then this follows from the fact that every
$a\in K^\dagger$ is a square. If $R$ is a CM-order, then
this follows from the fact that $R=R^+[F]$ is a free
$R^+$-module of rank two, so that $R/aR\cong(R^+/aR^+)^2$ as $R^+$-modules.

Finally, let $B(R)$ be the subgroup $\{P+\Pb : P\in G(\Mod_R)\}$ 
of $Z(R)$. We define $\CB(R)$ to be the group 
$Z(R)/(B(R)\cdot\Prin_R(K^\dagger))$.
Note that if $i\colon R\ra S$ is a homomorphism of
proper real/CM-orders then we get a group homomorphism
$i^*\colon\CB(S)\ra\CB(R)$ from the map $G(\Mod_S)\ra G(\Mod_R)$
obtained by viewing every $S$-module as an $R$-module. In fact,
$\CB$ is a contravariant functor from the category of real/CM-orders
to the category of abelian groups, and we will show in \S\SECgroup\ 
that $\CB(R)$ is a finite two-torsion group.

\proclaim{(\TTmain) Theorem}
Let $\CC$ be an isogeny class
of abelian varieties over a finite field and let $R=R_\CC$.
There is an element $I_\CC$ of $\CB(R)$
such that the elements of $G(\Mod_R)$ that are attainable in $\CC$ are 
precisely the effective elements of $Z(R)$ that map to $I_\CC$ in $\CB(R)$.
In particular, the isogeny class $\CC$ contains a 
principally polarized variety if and only if $I_\CC=0$.
\endproclaim

\demo{Proof} It is clear that attainable elements must be effective,
so the proof of Theorem~\TTmain\ splits up into three steps. 
In the first step we will show that every attainable element
of $G(\Mod_R)$ is contained in $Z(R)$.  In the second step
we will show that if $\lambda$ and $\mu$ are both polarizations
of varieties in $\CC$ then $\eps([\ker\lambda]_\CC)$ and
$\eps([\ker\mu]_\CC)$ differ by an element of 
$B(R)\cdot\Prin_R(K^\dagger)$.
In the third step we will show that if $P$ is an effective element of
$Z(R)$ that differs from an attainable element by an element
of $B(R)\cdot\Prin_R(K^\dagger)$, then $P$ is itself attainable.
The theorem will follow if we then take $I_\CC$ to be the image
in $\CB(R)$ of any attainable element of $Z(R)$.

{\it Step one.}
 Suppose $\lambda$ is a polarization of a variety
in $\CC$ and let $X$ be its kernel.
There is a non-degenerate alternating pairing $X\times X\ra \Gm$
(see \cite{\refMu\rm, \S23}), so there
is an isomorphism between $X$ and its Cartier dual $\Xh$.
Let $\Xll$ denote the local-local part
of $X$ and let $Y$ denote the non-local-local part of $X$.
Then $\Xll$ and $Y$ are both self-dual.  From \S\SECGrot\
we know that $\widehat{\CP(Y)}=\CP(\Yh)$
(where $\CP$ is the functor defined in \S\SECGrot)
and it follows easily that 
$\eps([X]_\CC)=\eps([\Xh]_\CC)=\ovl{\eps([X]_\CC)}$.
Therefore, attainable elements are symmetric.

Let $M$ be the $R$-module $\CP(Y)$.  The non-degenerate alternating
pairing from $Y$ to $\Gm$ gives us a non-degenerate alternating
$R$-semi-balanced
pairing from $M$ to $\BQ/\BZ$, and we can view this pairing as 
an $R^+$-balanced pairing on the $R^+$-module $M$. Lemma~\LLbalanced\ below
shows that every simple $R^+$-module occurs in $M$ an even number of times,
so $\eps([Y]_\CC)$ is in the kernel of the map 
$G(\Mod_R)\ra G(\Mod_{R^+})\otimes(\BZ/2\BZ)$.
We know from Riemann-Roch (see \cite{\refMu\rm, \S16})
that the rank of $X$ is a square,
and from what we have just seen the rank of $Y$ is a square. Therefore
the rank of $\Xll$ is a square, which shows that the simple group scheme
$\Bap$ occurs an even number of times in $\Xll$.  
It follows that $\eps([\Xll]_\CC)$ is in the kernel of 
$G(\Mod_R)\ra G(\Mod_{R+})\otimes(\BZ/2\BZ)$,
so $\eps([X]_\CC)$ is in this kernel as well.
Thus step one will be completed if we prove the following lemma.

\proclaim{(\LLbalanced) Lemma}
Let $A$ be a commutative ring and $M$ an 
$A$-module with $\#M<\infty$.  Suppose there is a non-degenerate balanced
alternating pairing $e\colon M\times M\ra\BQ/\BZ$.  
Then every simple $A$-module
that occurs in the Jordan-H\"older decomposition of $M$ occurs
an even number of times.
\endproclaim

\demo{Proof}
We will prove the lemma by induction on the cardinality of $M$.
The lemma is certainly true if $M$ is the zero module.
So suppose the lemma is true for all modules with cardinality less
than that of $M$.  We will show that it holds for $M$ also.

Let $S$ be a simple $A$-module that occurs in $M$, let $\pid$ be the
maximal ideal of $A$ that annihilates $S$, and let $L$ be the finite
field $A/\pid$. Identify $S$ with any
one-dimensional $L$-subspace of the $\pid$-torsion of $M$.
Then $e$ gives us a balanced alternating pairing $S\times S\ra \BQ/\BZ$ on
the $L$-vector space $S$. An easy lemma shows that there are no
non-degenerate balanced alternating pairings on a one-dimensional
vector space over a finite field
(see the proof of \cite{\refH\rm, Lem.~7.3, p.~2378}), 
so $e$ restricted to $S$ must be degenerate
and hence identically zero.  In other words, $S$ is an isotropic submodule
of the $A$-module $M$.

Let $S'\subset M$ be the annihilator of $S$ under $e$, so that
$S\subset S'$ by what we have just shown.  
The non-degenerate balanced pairing
$S\times M/S'\ra\BQ/\BZ$ induced from $e$ shows that $M/S'$ 
and $S$ are both simple modules killed by $\pid$, so they are isomorphic.
Thus in the Grothendieck group
of $A$ we have $[M]_A = [S'/S]_A + 2[S]_A$.
The pairing $e$ restricts to a non-degenerate balanced alternating pairing
$S'/S \times S'/S \ra\BQ/\BZ$.  The induction hypothesis shows
that every simple module that occurs in $S'/S$ occurs an even number
of times.  Therefore, the lemma holds for $M$, and we are done.
\myqed\enddemo

The center $K_\CC$ of the endomorphism algebra $E_\CC$ is a product
of number fields, say $K_\CC=K_1\times\cdots\times K_r$,
and we may write $E_\CC=E_1\times\cdots\times E_r$, where
each $E_i$ is a central simple $K_i$-algebra.  If $\lambda$
is a polarization of a variety in $\CC$, then the Rosati
involution~$\rosati$ on $E_\CC$
associated to $\lambda$ splits into an involution on each of the $E_i$.
If $K_i$ is a CM-field then $\rosati$ is an involution of the second
kind on $E_i$.  If $K_i$ is totally real, then $E_i$ ramifies at all
real primes of $K_i$ and $\rosati$ is a positive involution on $E_i$.
We will use these facts in steps two and three.

{\it Step two.}
Suppose $\lambda\colon A\ra\Ah$ and $\mu\colon B\ra \Bh$ are polarizations
of varieties in $\CC$ and let $\varphi\colon A\ra B$ be an isogeny
from $A$ to $B$.  Let $\nu$ be the polarization $\hat{\varphi}\mu\varphi$
of $A$, where $\hat{\varphi}\colon\Bh\ra\Ah$ is the
dual isogeny of $\varphi$.
Let $n$ be any positive integer such that $\ker\lambda\subset\ker(n\nu)$
as group schemes. Then there is an isogeny $\alpha\colon\Ah\ra\Ah$
such that $n\nu=\alpha\lambda$.  
A polarization is equal to its own dual isogeny, so we can equate
the right-hand side of the last equality with its dual to get
$n\nu=\lambda\hat{\alpha}$.
Using \cite{\refMu\rm, \S21, Application~III, pp.~208--210} (see
especially the final paragraph)
and the fact that $n\nu$ and $\lambda$ are polarizations,
we find that $\hat{\alpha}\in \End A$
is fixed by the Rosati involution and is totally positive.

The equality $n\hat{\varphi}\mu\varphi=\lambda\hat{\alpha}$ translates
into the equality
$$[\ker n]_\CC+[\ker\hat{\varphi}]_\CC+[\ker\mu]_\CC+
      [\ker\varphi]_\CC = [\ker\lambda]_\CC+[\ker\hat{\alpha}]_\CC$$
in $G(\Ker_\CC)$.  Let $Q=\eps([\ker\varphi]_\CC)$, so that
$\Qb=\eps([\ker\hat{\varphi}]_\CC)$.  Applying the isomorphism $\eps$ to the
equality above and using Theorem~\TTGrotnorms, we find that
$$\Prin_R(\Nrd(n)) + \Qb + \eps([\ker\mu]_\CC) + Q 
     = \eps([\ker\lambda]_\CC) + \Prin_R(\Nrd(\hat{\alpha})).$$
Now $n$ and $\hat{\alpha}$ are both fixed by the Rosati involution
and are totally positive, so by Theorems~\TTNrd\ and~\TTNrdfirst\ 
their reduced norms lie in the subgroup $K^\dagger$ of $K^*$.
Since $Q+\Qb$ is an element of $B(R)$, it is clear that
$\eps([\ker\lambda]_\CC)$ and $\eps([\ker\mu]_\CC)$ differ by 
an element of $B(R)\cdot\Prin_R(K^\dagger)$.

{\it Step three.}
Now suppose $P$ is an effective element of $Z(R)$ such that 
$$ P+Q+\Qb=\eps([\ker\lambda]_\CC) + \Prin_R(a) $$
for some $Q\in G(\Mod_R)$, some $a\in K^\dagger$, and some
polarization $\lambda\colon A\ra\Ah$ of a variety in~$\CC$.
By Theorems~\TTNrd\ and~\TTNrdfirst,
there is an $\alpha\in(\End A)\otimes\BQ$
that is fixed by the Rosati involution associated to $\lambda$, that
is totally positive, and such that $\Nrd(\alpha)=a$.
If we replace $\alpha$ by an integer multiple of itself and change $Q$
accordingly, we can assume that we have
$$ P+Q+\Qb=\eps([\ker\lambda]_\CC) + \Prin_R(\Nrd(\alpha)) $$
where $Q$ is effective and where $\alpha$ is a 
totally positive isogeny of $A$ that is fixed by the Rosati involution.
Then \cite{\refMu\rm, \S21, Application~III, pp.~208--210} tells
us that $\nu=\lambda\alpha$ is also a polarization of $A$,
and we have 
$$P+Q+\Qb=\eps([\ker\nu]_\CC). \eqno{\EEkernels}$$

Let $X=\ker\nu$ and let $e\colon X\times X\ra \Gm$ be the non-degenerate
alternating pairing on $X$ whose existence is shown in \cite{\refMu\rm, \S23}.
Let $Y$ be an element of $\Ker_\CC$ such that $S=\eps([Y]_\CC)$
is a simple $R$-module that occurs in the effective element $Q$ of 
$G(\Mod_R)$.
Suppose we can find an embedding of $Y$ into $X$ such that the pairing
$e$ restricted to $Y\times Y$ is the trivial pairing,
and let $\varphi$ be the natural isogeny from $A$ to $B=A/Y$.  
Then \cite{\refMu\rm, \S23, Cor.~to Thm.~2, p.~231} shows that there
is a polarization $\nu'$ of $B$ such that $\nu=\hat{\varphi}\nu'\varphi$.
In $G(\Mod_R)$ this gives us the equality 
$$\eps([\ker\nu]_\CC)= S + \eps([\ker\nu']_\CC) + \Sbar.$$
If we replace $Q$ by $Q-S$ and $\nu$ by $\nu'$, we will again
have equality~\EEkernels, but we will have decreased the number of
simple $R$-modules that occur in $Q$.  By applying this argument
repeatedly, we can finally reduce equality~\EEkernels\ to the
desired equality $P=\eps([\ker\nu]_\CC)$ for a polarization $\nu$
of a variety in $\CC$.  Thus step three and the proof of Theorem~\TTmain\ 
are completed by the following proposition. 
\myqed\enddemo

\proclaim{(\PPisotropic) Proposition}
Suppose $X$ and $Y$ 
are finite commutative group schemes over a finite field $k$, and suppose
$Y$ is simple.  Suppose further that there is a non-degenerate
alternating pairing $e\colon X\times X\ra\Gm$ and that $[X]-[Y]-[\Yh]$
is an effective element of the Grothendieck group of finite
commutative group schemes over $k$.  Then there is an embedding of
$Y$ into $X$ such that $e$ restricted to $Y\times Y$ is the trivial
pairing.
\endproclaim

\demo{Proof}
The pairing $e$ splits into three non-degenerate alternating
pairings, one each on the reduced-reduced part of $X$,
the product of the reduced-local and the local-reduced parts of $X$,
and the local-local part of $X$. 
Therefore it is enough to prove the theorem when $X$ is either a
reduced-reduced group, a product of a reduced-local and a local-reduced 
group, or a local-local group.  Clearly we may assume that $Y$ is of the
same type as $X$.

Suppose $X$ is reduced-reduced.  Let $S$ be the group ring over $\BZ$
of the Galois group $\Gal(\kb/k)$; the map $\sigma\mapsto\sigma^{-1}$
of the Galois group gives us an involution $\cxcj$ on $S$, and
$S$ is commutative since $k$ is finite.
The group schemes $X$ and $Y$ are completely determined by the
$S$-modules $X(\kb)$ and $Y(\kb)$, and the pairing on $X$ gives 
us a pairing $X(\kb)\times X(\kb)\ra \BQ/\BZ$ that is non-degenerate,
alternating, 
and semi-balanced with respect to the involution $\cxcj$ on~$S$.
The proposition follows in this case from \cite{\refH\rm, Prop.~7.1, p.~2378}.

Suppose $X$ is a product of a reduced-local and a local-reduced
group scheme.  We may assume that $Y$ is reduced-local.
Let $i\colon Y\ra X$ be any embedding of $Y$ into $X$; it is not hard to
show that such embeddings exist.  Then $e$ restricted to $Y\times Y$ 
cannot be non-degenerate, because otherwise we would have an isomorphism
between the reduced-local group scheme $Y$ and its local-reduced dual $\Yh$.
Therefore $e$ is trivial on $Y\times Y$.

Suppose $X$ is local-local, so that $Y$ must be isomorphic to $\Bap$,
where $p=\charact k$.
If $p>2$ then a straightforward
Hopf algebra computation shows that
there are no non-degenerate alternating pairings on $\Bap$, so
again we can take any embedding of $\Bap$ into $X$. 
We are left with the case where $p=2$ and $Y=\alphatwo$.

Let $F=\End\alphatwo\cong k$, let $U$ be the $F$-vector space
$\Hom(\alphatwo,X)$ (where $F$ acts by premultiplication), and
let $V$ be the $F$-vector space of alternating pairings on
$\alphatwo$ (where $F$ acts by premultiplication on the first
factor). A computation shows that $V$ is one-dimensional.
Let $i\colon \alphatwo\ra X$ be any embedding of $\alphatwo$ into $X$
and let $j$ be any embedding of $\alphatwo$ into the
kernel of the composite map 
$X\cong\Xh\ra\widehat{\alphatwo}$, 
where the isomorphism $X\cong\Xh$ is obtained
from the pairing $e$ and where $\Xh\ra\widehat{\alphatwo}$ is dual to $i$. 
If $j$ is a multiple of $i$ in $U$ then
the pairing $e$ restricted to $i(Y)\times i(Y)$ is trivial and
we are done.  
Otherwise, let $p_1$ and $p_2$ be the elements $e\circ(i\times i)$
and $e\circ(j\times j)$ of $V$.  If either $p_1$ or $p_2$ is trivial
we are done, so assume  they are not.  Using the fact that $i(\alphatwo)$
and $j(\alphatwo)$ are orthogonal under $e$, one can show that
for every $\beta\in F$ we have 
$$e\circ((i+\beta j)\times(i+\beta j))=p_1+\beta^2 p_2.$$
Since $p_1$ and $p_2$ are non-zero elements of the one-dimensional
$F$-vector space $V$
and since the squaring map $F\ra F$ is onto, we find that there is
a choice of $\beta$ so that the pairing $e$ restricted to the
image of the embedding $i+\beta j$ is trivial.
\myqed\enddemo

\head \SECgroup.  The obstruction group
\endhead 
In this section we will determine how one can calculate the
obstruction group $\CB(R_\CC)$ for an isogeny class $\CC$ 
of abelian varieties over a finite field $k$.  
We begin by calculating $\CB(R)$ in the easiest cases: when $R$ is 
totally real, and when $R$ is the ring of integers of a CM-field.

\proclaim{(\PPreal) Proposition}
Suppose $R$ is an order in a 
product of totally real number fields.  Then $\CB(R)\cong 0$.
\endproclaim

\demo{Proof} In this case the involution on $R$ is trivial 
and $R=R^+$.
We check that $Z(R)=2G(\Mod_R)$ and $B(R)=2G(\Mod_R)$.
Then clearly $\CB(R)=0$. 
\myqed\enddemo

\proclaim{(\PPCM) Proposition}
Let $\Oid$ be the ring of integers of a CM-field $K$.
Then $$\CB(\Oid)\cong \Cl^+(K^+)/N_{K/K^+}(\Cl(K)),$$
where $\Cl^+(K^+)$ denotes the narrow class group of $K^+$.
If $K/K^+$ is ramified at a finite prime, then 
$\Cl^+(K^+)/N_{K/K^+}(\Cl(K))$ is the zero group. Otherwise,
the Artin map provides an isomorphism
$$\Cl^+(K^+)/N_{K/K^+}(\Cl(K))\cong\Gal(K/K^+).$$
\endproclaim

\demo{Proof} 
The map $\aid\mapsto\Oid/\aid\Oid$ gives an isomorphism between the
ideal group of $K^+$ and $Z(\Oid)$.  Under this isomorphism, the subgroup
of principal ideals of $K^+$ that are generated by totally positive elements
is identified with $\Prin_\Oid(K^\dagger)$.  Also, the subgroup
of norms of ideals of $K$ is identified with $B(\Oid)$.
Therefore $$\Cl^+(K^+)/N_{K/K^+}(\Cl(K))\cong\CB(\Oid).$$
The rest of the proposition is \cite{\refH\rm, Prop.~10.1, p.~2385}. 
\myqed\enddemo

\proclaim{(\CCodddegree) Corollary}
Let $K$ be a CM-field.
If $[K^+:\BQ]$ is odd then $\CB(\Oid_K)\cong 0$.
\endproclaim

\demo{Proof} We see from \cite{\refH\rm, Lem.~10.2, p.~2385}
that if $[K^+:\BQ]$ is odd then 
$K/K^+$ must be ramified at a finite prime.  The corollary follows
from Proposition~\PPCM.
\myqed\enddemo

The next proposition allows us to compute the obstruction group for
one ring in terms of that of an overring.

\proclaim{(\PPpushout) Proposition}
Suppose $R$ and $S$ are proper real/CM-orders and $R$ 
is a subring of finite index in $S$.  Then
$$\matrix
    Z(S)/B(S)     & \lra  & \CB(S)           \\
\Btbmapdown{N}    &       & \Btbmapdown{i^*} \\
    Z(R)/B(R)     & \lra  & \CB(R)           \\
\endmatrix$$
is a push-out diagram, where $i\colon R\ra S$ is the inclusion map,
where $N$ is induced from the norm from $G(\Mod_S)$ to $G(\Mod_R)$,
and where the horizontal maps are the natural reduction maps.
\endproclaim

\demo{Proof}
We start with the exact sequence
$$Z(S)\lllongmapright{(N,-1)} Z(R)\oplus Z(S)\lllongmapright{1\oplus N} 
                          Z(R)\lllra 0$$
where $N$ is the norm map from $Z(S)$ to $Z(R)$. Let $K=R\otimes\BQ$;
since $R$ is of finite index in $S$ we have $K=S\otimes\BQ$. From 
\cite{\refH\rm, Lem.~2.1, p.~2365} we have $\Prin_R=N\circ \Prin_S$,
and dividing the last two terms of the preceding sequence
by the image of $K^\dagger$
under $\Prin_S$ and $\Prin_R$ we get the exact sequence
$$Z(S)\lllra Z(R)\oplus (Z(S)/\Prin_S(K^\dagger))
    \lllra Z(R)/\Prin_R(K^\dagger)\lllra 0.\eqno{\EEbottom}$$
We also have an exact sequence
$$B(S)\lllongmapright{(N,-1)} B(R)\oplus B(S)\lllongmapright{1\oplus N} 
                          B(R)\lllra 0.\eqno{\EEtop}$$
If we take the cokernel of the natural map from \EEtop\ to \EEbottom\ we
get the exact sequence
$$Z(S)/B(S)\lllongmapright{(N,-\chi)} 
    (Z(R)/B(R))\oplus \CB(S) \lllongmapright{\psi\oplus i^*} \CB(R)
      \lllra 0,\eqno{\EEquot}$$
where $\chi\colon Z(S)/B(S)\ra \CB(S)$ and $\psi\colon Z(R)/B(R)\ra\CB(R)$
are the natural reduction maps.  This sequence being exact 
is precisely what it means for the diagram of the proposition to be
a push-out diagram.
\myqed\enddemo

The groups on the left-hand side of Proposition~\PPpushout\ can
be replaced by finite groups. We need some notation to do this. 
For every proper real/CM-order $R$ let $H(R)=Z(R)/B(R)$.
The group $H(R)$ is a vector space over the field with two elements, and 
a basis for $H(R)$ as an $\Ftwo$-vector space is given
by the images in $H(R)$ of the
classes of the simple $R$-modules of the form $R/\pid$, where $\pid$
is a maximal ideal of $R$ that is fixed by the involution of $R$
and whose residue field has even degree over the residue field
of the prime of $R^+$ that lies under it; we call such maximal
ideals {\it generating primes} of $H(R)$.

Now suppose $R$ and $S$ are as in Proposition~\PPpushout.  Let $X$ denote
the set of generating primes of $H(R)$, and for each $\pid$
in $X$ let $x_\pid$ be the image of the simple $R$-module $R/\pid$ 
in $H(R)$.  Similarly, let $Y$ denote the 
set of generating primes of $H(S)$, 
and for each $\qid$ in $S$ let $y_\qid$ be
the image of $S/\qid$ in $H(S)$.  Let $\Xgood$ be the set of
$\pid$ in $X$ such that $S\otimes_R R/\pid\cong R/\pid$, and
let $\Ygood$ be the set of primes $\qid$ of $Y$ such that
$S\otimes_R R/(R\cap\qid)\cong R/(R\cap\qid)$.
It is not hard to see that every prime $\pid\in\Xgood$ has exactly
one prime $\qid$ of $S$ lying over it, that this prime of $S$ 
lies in $\Ygood$, and that the residue fields of $\pid$ and $\qid$ are
equal. Likewise, it is easy to see that every prime $\qid\in\Ygood$
lies over a prime $\pid\in\Xgood$.
  
Let $\Cgood$ and $\Cbad$ be the subspaces of $H(R)$ spanned by the sets 
$\{x_\pid \ |\ \pid\in\Xgood\}$ and 
$\{ x_\pid \ |\ \pid\in X\setminus\Xgood\}$,
respectively, and let $\Dgood$ and $\Dbad$ be the subspaces of $H(S)$
spanned by the sets $\{y_\qid \ |\ \qid\in\Ygood\}$ and
$\{ y_\qid \ |\ \qid\in Y\setminus\Ygood\}$, respectively.
It is easy to see from what we have noted above that the norm map
from $H(S)$ to $H(R)$ induces an isomorphism between $\Dgood$ and
$\Cgood$, and that the image of $\Dbad$ lies in $\Cbad$. 

\proclaim{(\PPpushouttwo) Proposition}
Let notation and assumptions be
as in Proposition~\rom{\PPpushout}. Then
$$\matrix
     \Dbad         & \lra  & \CB(S)          \\
 \Btbmapdown{N}    &       & \Btbmapdown{i^*}\\
     \Cbad         & \lra  & \CB(R)          \\
\endmatrix$$
is a push-out diagram.
\endproclaim

\demo{Proof} We would like to show that the sequence
$$\Dbad\lllongmapright{(N,-\chi)} 
    \Cbad\oplus \CB(S) \lllongmapright{\psi\oplus i^*} \CB(R)
      \lllra 0\eqno{\EEfinite}$$
is exact.  Consider the natural map from this sequence to
the sequence~\EEquot.  We have already noted that the norm from 
$H(S)$ to $H(R)$ induces an isomorphism between $\Dgood$ and $\Cgood$.
A simple diagram chase making use of this fact shows that~\EEfinite\ is exact.
\myqed\enddemo

One can show that every $\pid\in X\setminus\Xgood$ is singular
and that every $\qid\in Y\setminus\Ygood$ lies over a singular
prime of $R$, so $\Dbad$ and $\Cbad$ are finite.  This leads
to the following corollary:

\proclaim{(\CCfinite) Corollary}
Let $R$ be a proper real/CM-order.
Then $\CB(R)$ is a finite two-torsion group.
\endproclaim

\demo{Proof} Let $K=R\otimes\BQ$ and write $K$ as a product of
fields $K_1\times\cdots\times K_r$.  For each $i$ let $\Oid_i$ be the
ring of integers of $K_i$, and let $S$ be the product of the $\Oid_i$.
It is clear that the functor $\CB$ applied to a product of rings
gives the sum of $\CB$ applied to each factor, so using 
Propositions~\PPreal\ and~\PPCM\ we can see that
$\CB(S)=\CB(\Oid_1)\oplus\cdots\oplus\CB(\Oid_r)$ is a finite two-torsion
group.  Now apply Proposition~\PPpushouttwo\ to $R$ and $S$.
Since the upper right and lower left groups in the diagram of 
Proposition~\PPpushouttwo\ are finite two-torsion groups, so is $\CB(R)$. 
\myqed\enddemo

We end this section with two comments.
Suppose $R=R_\CC$ for an isogeny class~$\CC$,
and let $R_{\roman{CM}}$ be the image of $R$ under the projection
map from $R\otimes\BQ$ to the product of the CM-fields occurring
in $R\otimes\BQ$.
One can use Proposition~\PPpushouttwo\ to show that the natural
map $\CB(R_{\roman{CM}})\ra\CB(R)$ is an isomorphism.
Also, the ring $R_{\roman{CM}}$
is itself of the form $R_{\CC'}$ for an isogeny class $\CC'$. 
In \cite{\refH\rm, \S9} it is shown how 
the CM-order $R_{\roman{CM}}$ is determined 
by the characteristic polynomial $h_{\CC'}$.
(It is assumed in~\cite{\refH}
that the isogeny classes contain only ordinary abelian varieties,
but the only way this assumption is used in 
\cite{\refH\rm, \S9} is in the fact
that the rings associated with ordinary isogeny classes are automatically
CM-orders.)
Thus, one can actually calculate the obstruction group for an isogeny
class $\CC$ if one knows $h_\CC$.
Finally, we note that although the definition of $\CB(R)$ found
in \cite{\refH} is different from the one given in this paper,
the two definitions agree when $R$ is a proper CM-order that is locally
free of rank two over $R^+$.  In particular, one can show that
the definitions agree when $R$ is a CM-order of the form $R_\CC$.

\head \SECelement.  Restricting the obstruction element
\endhead
Let $\CC$ be an isogeny class of
abelian varieties over a finite field. 
In this section we show that $I_\CC$ must lie
in a particular subgroup of $\CB(R_\CC)$.

\proclaim{(\PPOimage) Proposition}
Let $\Oid$ be the integral
closure of $\BZ$ in the product of fields $K_\CC$
and let $i\colon R_\CC\ra\Oid$ be the inclusion map.  
Then the obstruction element $I_\CC\in \CB(R_\CC)$ lies in the image of 
$\CB(\Oid)$ under the map $i^*\colon\CB(\Oid)\ra\CB(R_\CC)$.
\endproclaim

\demo{Proof}
As usual we let $R=R_\CC$. By \cite{\refW\rm, Thm.~3.13, p.~534},
there is a variety $A$ in $\CC$ such that $\Oid\subset\End(A)$.
Let $\lambda$ be any polarization of $A$ and let $X=\ker\lambda$.
We write $\Xll$ for the local-local part of $X$ and $Y$ for the
non-local-local part of $X$.  
Theorem~\TTmain\ shows that $I_\CC$ is the image of $X=\Xll\times Y$ 
in $\CB(R)$.

We showed in step one of the proof of Theorem~\TTmain\ that 
the image of $\Xll$
in $G(\Ker_\CC)$ is an even multiple of the simple group scheme $\Bap$.
Therefore the image of $\Xll$ in $G(\Mod_R)$ is an even multiple of
the simple $R$-module $\BZ/p\BZ$ on which $F$ and $V$ act as zero.
Thus the image of $\Xll$ in $G(\Mod_R)$ is in $B(R)$, so the image
of $\Xll$ in $\CB(R)$ is zero.  We are left to consider $Y$.

Let $n$ be any positive integer such that $R\supset n\Oid$ and let $\mu$
be the polarization $n\lambda$ of $A$. 
Let $\Oid$ act on $\Ah$ by sending
$a\in\Oid\subset\End A$ to the endomorphism $\widehat{a'}$ of $\Ah$, where
$\rosati$ is the Rosati involution of $\End A$ associated to $\lambda$.  
Then the definition of the Rosati involution shows that $\lambda$ and $\mu$
are $\Oid$-equivariant, so $\Oid$ acts on $X$ and on the kernel $X'$ of $\mu$.
Let $e$ and $e'$ be the non-degenerate alternating pairings
from $X\times X$ and $X'\times X'$ to $\Gm$ defined in \cite{\refMu\rm, \S23}.
The pairings $e$ and $e'$ are defined over $k$,
so for every $\alpha\in R$, for every $k$-algebra $S$, and for every
pair of $S$-valued points $x$ and $y$ of $X'$ we have 
$e'(\alpha x, y)=e'(x,\ovl{\alpha}y)$; 
we express this fact by saying that the pairing $e'$ is $R$-semi-balanced.
We also know that for every $k$-algebra $S$, if $x\in X(S)$ and
$y\in X'(S)$ then $e'(x,y)=e(x,ny)$; this is \cite{\refMu\rm, item~4, p.~228}.
Using this relation between the two pairings and the fact
that $e'$ is $R$-semi-balanced, one can show that $e$ is 
$\Oid$-semi-balanced; that is, for every $\alpha\in\Oid$,
for every $k$-algebra $S$, and for every pair of points
$x, y\in X(S)$, we have $e(\alpha x, y)=e(x, \ovl{\alpha} y)$.

The ring $\Oid$ acts on $Y$, and since there are no non-zero morphisms
between $Y$ and $\Xll$, the alternating $\Oid$-semi-balanced pairing
$e$ restricted to $Y\times Y$ is still non-degenerate.
Let $M=\CP(Y)$, where $\CP$ is the functor from $\Ker_\CC$ to
$\Mod_R$ that was defined in \S\SECGrot.  Since $\Oid$
acts on $Y$, the $R$-module $M$ is actually an $\Oid$-module, and 
our pairing from $Y\times Y$ 
to $\Gm$ gives us a non-degenerate alternating 
$\Oid$-semi-balanced pairing from $M\times M$ to $\BQ/\BZ$.
As in the similar situation in step one of the proof of Theorem~\TTmain,
this tells us that $[M]_\Oid$ is a symmetric element
of $G(\Mod_\Oid)$, and Lemma~\LLbalanced\ tells us that the image of
$M$ in $G(\Mod_{\Oid^+})\otimes(\BZ/2\BZ)$ is zero.  Thus
$[M]_\Oid$ is in $Z(\Oid)$. The natural map from $Z(\Oid)$
to $\CB(R)$ factors through $\CB(\Oid)$, so the image of $[M]_R$ in 
$\CB(R)$, which is the image of $[M]_\Oid$ in $\CB(R)$, lies in the
image of $\CB(\Oid)$.   This proves the proposition. 
\myqed\enddemo

The following proposition, which should be compared to
\cite{\refH\rm, Prop.~11.3, p.~2390}, gives sufficient conditions for
$i^*$ to be the zero map in the case where $\CC$ is simple.

\proclaim{(\PPzeroconditions) Proposition}
Let $\CC$ be an isogeny class of simple $g$-dimensional abelian varieties
over a finite field $k$, let $K=K_\CC$, let $\Oid$
be the ring of integers of $K$, and let $i\colon R_\CC\ra\Oid$
be the inclusion map. If $K$ is totally real then $i^*$ is the
zero map.  Suppose $K$ is a CM-field.
If $K/K^+$ is ramified at a finite prime,
or if there is a prime of $K^+$ that divides $(F-V)$ and that is inert
in $K/K^+$, then $i^*$ is the zero map.
Otherwise, $i^*$ is not the zero map and $g$ is even.
\endproclaim

\demo{Proof} 
If $K$ is totally real then $\CB(R)\cong 0$ by Proposition~\PPreal,
and there is nothing more to say.  So suppose $K$ is a CM-field.
If $K/K^+$ is ramified at a finite prime then $\CB(\Oid)\cong 0$ 
by Proposition~\PPCM\ and $i^*$ must be the zero map.  
We are left with the case where $K/K^+$ is
unramified at all finite primes; then 
Proposition~\PPCM\ shows that $\CB(\Oid)$
has two elements, so $i^*$ is the zero map precisely when it kills
the non-zero element of $\CB(\Oid)$.
We will use Proposition~\PPpushouttwo\ to determine
when this is the case.

Let us adopt the notation of Proposition~\PPpushouttwo\ and apply
the proposition with $S=\Oid$.  Let $I$ be the non-zero element of
$\CB(\Oid)$.
The generating primes of $H(\Oid)$ are the primes of $\Oid$ 
that are inert in $K/K^+$, and Proposition~\PPCM\ tells us that
for every prime $\qid$ of $\Oid$ that is inert in $K/K^+$,
the image of $y_\qid$ in $\CB(\Oid)$ is $I$.
Now, Proposition~\PPpushouttwo\ shows that $i^*(I)=0$ if and only if
there is a $y\in \Dbad$ that maps to $I$ in $\CB(\Oid)$ and
to zero in $\Cbad$; for the moment
let us call such a $y$ an {\it annihilating} element.
Since the basis elements $y_\qid$ of $\Dbad$
all map to $I$ in $\CB(\Oid)$, and since the image of $y_\qid$
in $\Cbad$ is either $0$ or a basis element $x_\pid$, we see
that there is an annihilating $y\in\Dbad$ if and only if there
is an annihilating $y$ of the form $y=y_\qid$.  
One can check that an element $y_\qid$ maps
to zero in $\Cbad$ if and only if the residue field of $\qid$
has even degree over the residue field of the prime $\pid=\qid\cap R$
of $R$, and this last condition is equivalent to the
condition that complex conjugation on $\Oid/\qid$ act trivially
on $R/\pid\subset \Oid/\qid$.  Since $R=\BZ[F,V]$, this will
be the case if and only if $F$ and its complex conjugate $V$ are
congruent to one another modulo $\qid$, 
if and only if $\qid$ divides $F-V$.
Thus, $i^*$ is the zero map if and only if there is a prime
of $K^+$ that is inert in $K/K^+$ and that divides $F-V$.
Finally, if $i^*$ is not the zero map then Corollary~\CCodddegree\ shows
that $[K^+:\BQ]$ is even, and since $g$ is a multiple of this degree, 
$g$ is even also. 
\myqed\enddemo

We close this section by noting that Theorems~\TTmainPP\ and~\TTodddim\ 
follow immediately from Proposition~\PPzeroconditions,
Proposition~\PPOimage, and Theorem~\TTmain.

\enddocument

\Refs

\ref\no\refD
\by       P. Deligne
\paper    Vari\'et\'es ab\'eliennes ordinaires  sur un corps fini
\jour     Invent. Math. 
\vol      8
\yr       1969
\pages    238--243
\endref

\ref\no\refH
\by       E. W. Howe
\paper    Principally polarized ordinary abelian varieties over finite fields
\jour     Trans. Amer. Math. Soc.
\vol      347
\yr       1995
\pages    2361--2401
\endref

\ref\no\refMu
\by       D. Mumford
\book     Abelian Varieties
\bookinfo Tata Inst. Fund. Res. Stud. Math.~{\bf 5}
\publ     Oxford University Press
\publaddr Oxford
\yr       1985
\endref

\ref\no\refO
\by       T. Oda
\paper    The first de Rham cohomology group and Dieudonn\'e modules
\jour     Ann. Sci. \'Ecole Norm. Sup. (4) 
\vol      2
\yr       1969
\pages    63--135
\endref

\ref\no\refPR
\by       V. Platonov and A. Rapinchuk
\book     Algebraic Groups and Number Theory
\bookinfo Pure Appl. Math. {\bf 139}
\publ     Academic Press
\publaddr San Diego, California
\yr       1994
\endref

\ref\no\refRam
\by       C. P. Ramanujam
\paper    The theorem of Tate
\paperinfo appendix 1 of \cite{\refMu}
\endref

\ref\no\refR
\by       I. Reiner
\book     Maximal Orders
\bookinfo London Math. Soc. Monographs {\bf 5}
\publ     Academic Press
\publaddr London
\yr       1975
\endref

\ref\no\refSch
\by       W. Scharlau
\book     Quadratic and Hermitian Forms
\bookinfo Grundlehren Math. Wiss. {\bf 270}
\publ     Springer-Verlag
\publaddr Berlin
\yr       1985
\endref

\ref\no\refSe
\by       J.-P. Serre
\book     Local Fields
\bookinfo Grad. Texts in Math.~{\bf 67}
\publ     Springer-Verlag
\publaddr New York
\yr       1979
\endref

\ref\no\refT
\by       J. Tate
\paper    Classes d'isog\'enie des vari\'et\'es ab\'eliennes 
          sur un corps fini (d'apr\`es T. Honda)
\paperinfo expos\'e {\bf 352}
\inbook   S\'eminaire Bourbaki 1968/69
\bookinfo Lecture Notes in Math. {\bf 179}
\publ     Springer-Verlag
\publaddr Berlin
\yr       1971
\pages    95--110
\endref

\ref\no\refW
\by       W. C. Waterhouse
\paper    Abelian varieties over finite fields
\jour     Ann. Sci. \'Ecole Norm. Sup. (4)
\vol      2
\yr       1969
\pages    521--560
\endref


\endRefs
\bye